\pdfoutput=1
\documentclass[pdflatex,sn-basic]{sn-jnl}

\usepackage{lmodern}%
\usepackage{graphicx}%
\usepackage{multirow}%
\usepackage{amsmath,amssymb,amsfonts}%
\usepackage{amsthm}%
\usepackage{mathrsfs}%
\usepackage{mathtools}%
\usepackage[title]{appendix}%
\usepackage{xcolor}%
\usepackage{textcomp}%
\usepackage{manyfoot}%
\usepackage{booktabs}%

\begin{document}

\title[Exit law of a planar OU process on annular sectors]{The joint exit-time and exit-location law of a planar Ornstein--Uhlenbeck process on annular sectors and principal-axis rectangles}

\author[1]{\fnm{Tristan} \sur{Guillaume}}\email{tristan.guillaume@cyu.fr}

\affil[1]{\orgdiv{Laboratoire Thema}, \orgname{CY Cergy Paris Universit\'e}, \orgaddress{\street{33 boulevard du port}, \city{Cergy-Pontoise Cedex}, \postcode{F-95011}, \country{France}}}

\abstract{We study the first exit of an isotropic planar Ornstein--Uhlenbeck process from an annular sector, the region bounded by two concentric circular arcs and two radial segments, and obtain explicit eigenfunction expansions for the associated exit functionals. After the ground-state transformation that renders the generator self-adjoint, polar separation reduces the problem to a radial confluent hypergeometric equation whose two independent (Whittaker) solutions of generally non-integral order are combined through a two-radius determinant that fixes the spectrum. In this way we obtain the survival probability, the density and moments of the exit time, and --- the principal contribution --- the joint law of the exit time and the exit boundary, which resolves both the instant of exit and the boundary piece through which it occurs. As a companion separable case, we also treat a genuinely correlated, reversible planar Ornstein--Uhlenbeck process on a principal-axis rectangle, where the radial functions are replaced by parabolic cylinder functions: the survival probability factorizes into one-dimensional problems, while the joint law of the exit time and the exit side does not. The expansions require only one-dimensional root-finding and quadrature and are validated against Monte Carlo simulation. An application to an optically trapped colloidal particle is discussed, a setting in which the annular geometry arises naturally.}



\maketitle

\section{Introduction}\label{sec:1}

The computation of the distribution of the first exit time of a diffusion process from a bounded domain is a recurring problem throughout the mathematical sciences, arising in fields as varied as physics, finance, neuroscience, reliability engineering and population dynamics. When the underlying process is a Brownian motion, a substantial body of exact and approximate results is available, for one-sided as well as two-sided barriers, and for constant as well as moving boundaries, including stochastic ones. The situation is markedly different for mean-reverting processes. The Ornstein--Uhlenbeck (OU) process --- after Brownian motion the most widely used diffusion in applied modeling, and, by Doob's theorem, the only stationary Gauss--Markov process --- has proved considerably more resistant to exact first-passage analysis. Even in one dimension, no elementary closed form is available for the density of the first passage time of an OU process to a fixed level, and the most explicit representations rest on spectral expansions involving parabolic cylinder functions and the zeros thereof \citep{alili2005,darling1953,linetsky2004}. The two-sided exit from an interval has been treated along the same spectral lines, together with the moments and the Laplace transform of the exit time \citep{dinardo2001}.

In higher dimension, exact results are scarcer still, and --- this is the point of departure of the present article --- those that exist are confined to domains whose boundary is everywhere of a single coordinate type. The most complete is due to \citet{grebenkov2015}, who obtains the mean exit time, the moment-generating function and the survival probability of a multidimensional OU process from a ball, expressed through confluent hypergeometric functions; the high-dimensional asymptotics of the mean exit time from a ball have since been examined in~\citet{kersting2023}. For a ball, rotational invariance reduces the problem to its radial part, and the angular variables play no role. At the other extreme, the first-exit problem for a planar OU process through a smooth closed curve, namely a time-varying ellipse, has been addressed by means of the Laplace transform of the exit-time density and its numerical inversion \citep{dicrescenzo2024}, while the inverse problem --- recovering the boundary that yields a prescribed exit-time law --- has been treated numerically, with applications to neuronal modeling \citep{civallero2019}. What these contributions share is that the boundary is everywhere of one kind: either purely radial, as for the ball, or a single smooth closed curve, as for the ellipse. To the best of our knowledge, no explicit, boundary-resolved result is available for an OU process on a domain whose boundary combines curved radial portions with rectilinear angular portions --- that is, on a domain possessing corners.

The present article provides such a result. We consider an isotropic planar OU process and study its first exit from the annular sector

\begin{equation}
D = \{(r,\varphi):\ 0 < r_{1} < r < r_{2},\ {0 \leq \varphi}_{1} < \varphi < \varphi_{2} \leq 2\pi\},
\tag{1.1}
\end{equation}

bounded by two concentric circular arcs and two radial segments. The annular sector is, in a precise sense, the minimal planar domain that simultaneously carries both kinds of boundary while remaining clear of the coordinate singularity at the origin; the strict positivity of the inner radius \(r_{1}\) is, as will be seen, not a technical convenience but the source of the genuinely two-sided radial structure of the problem. The domain (1.1) moreover shares its radial Kummer equation with the disk (\(r_{1} \rightarrow 0\), \(\varphi_{2} - \varphi_{1} = 2\pi\)), the full annulus (\(\varphi_{2} - \varphi_{1} = 2\pi\)), the circular wedge (\(r_{1} \rightarrow 0\), \(r_{2} \rightarrow \infty\)) and the bounded sector. The circular wedge and bounded sector are genuine Dirichlet-sector limits, retaining the absorbing radial edges; the disk and full annulus, as explained in Remark 3, are recovered only once those edges are replaced by periodic angular conditions, the lowest angular order of the sector itself tending to 1/2 rather than to 0. With that proviso, a single formula on (1.1) unifies a family of exit problems that would otherwise be handled separately.

The central object of the paper is not merely the survival probability, but the joint law of the exit time and the exit boundary: the probability that the process leaves \(D\) after a prescribed time, together with the identification of which of the four boundary pieces is struck, and at what point. For a rotationally symmetric domain such as the ball, this refinement is vacuous --- there is a single boundary, and the question of which boundary is reached carries no information --- and it is precisely the mixed, cornered geometry of the annular sector that renders it meaningful, calling on the full angular as well as radial structure of the problem. More generally, we obtain a closed form for a discounted functional of the exit time and exit location, from which the survival probability, the exit-time density and the boundary-hitting probabilities all follow as particular cases.

The analysis proceeds along the classical connection between diffusions and parabolic partial differential equations \citep{karatzas1991,stroock1979}: the sought probability solves a backward Kolmogorov equation on \(D\) with a constant initial condition and homogeneous Dirichlet data on \(\partial D\). A ground-state (similarity) transformation removes the first-order drift and converts the generator into the Hamiltonian of a two-dimensional harmonic oscillator, after which separation of variables in polar coordinates yields three ordinary differential equations. The angular equation produces a Fourier sine series whose orders \(\nu_{m} = m\pi/\left( \varphi_{2} - \varphi_{1} \right)\) are generally non-integral. The radial equation is a confluent hypergeometric equation, whose two independent solutions are Whittaker functions of generally non-integral order; the radial eigenvalues are then the roots of a transcendental equation expressing the vanishing of a two-radius determinant of these functions. The coefficients of the expansion are recovered by a generalized Fourier projection of the constant initial datum against the eigenfunctions, taken with respect to the Gaussian weight under which the OU generator is self-adjoint. The resulting series is a double sum of products of elementary and Whittaker functions; as with the confluent-hypergeometric expansions of \citet{grebenkov2015} and those of the confined-oscillator and exactly-soluble quantum-ring literature \citep{lewyanvoon2003,tan1996}, the survival and density series converge rapidly for positive times, while the boundary-integrated quantities require the additional care usual for Dirichlet eigenfunction expansions near absorbing boundaries (Section 5). In all cases the terms are evaluated with standard scientific computing software, so that no recourse to finite-difference or other space-discretization schemes is required.

It must be emphasized that, for the circular annular-sector geometry, the polar separation underlying this construction requires the centered OU generator to be isotropic. As soon as the process is genuinely correlated --- as soon, that is, as its diffusion or drift matrix fails to be a scalar multiple of the identity --- the ground-state transformation produces an anisotropic harmonic potential, which separates along its principal Cartesian axes rather than in the polar coordinates that the circular boundary would require. Conversely, a correlated OU process does separate, but on a rectangle aligned with its principal axes; the pertinent special functions are then parabolic cylinder functions, and the decorrelation is effected by a change of variables. We therefore complement the annular-sector result with the solution of the exit problem for a correlated OU process from such a rectangle, thereby covering the two canonical separable geometries: curved and isotropic on the one hand, rectilinear and correlated on the other. Since the bare survival probability on the rectangle factorizes into one-dimensional exit problems, the genuinely two-dimensional content there resides in the joint exit-time-and-location law and in the discounted functionals, which do not factorize.

Beyond its intrinsic interest as a completion of the catalogue of solvable OU exit problems, the result is motivated by settings in which the annular sector is the natural event region rather than an imposed one. The overdamped motion of a micron-sized particle in a harmonic optical trap is an OU process, and the experimental geometries of single-particle tracking --- apertures, channels and shells --- are precisely those in which escape occurs through an angular gap or across an annular boundary rather than across a sphere \citep{grebenkov2015}. The isotropic planar OU process is likewise a canonical local model for the noisy fluctuations of a dynamical system about a stable equilibrium with isotropic linear restoring force; in the polar (amplitude--phase) coordinates natural to such a setting, the event that the radial amplitude remains within a tolerance band while the angular coordinate stays within a prescribed range is exactly an exit from an annular sector. A true noisy stable focus carries in addition a rotational drift component, which renders the generator non-self-adjoint and so falls outside the polar separation used here; it is a natural extension, recorded among the open problems of Section 7. In all of these settings, the availability of an explicit formula confers the customary advantages over a purely numerical treatment: sensitivities with respect to every parameter follow by direct differentiation rather than repeated simulation, and the explicit formula furnishes a benchmark against which numerical --- and, increasingly, neural-network-based --- solvers may be calibrated, a role for which exactly solvable problems on cornered geometries are in short supply and in which standard discretization schemes are least reliable, owing to the short-time leakage of probability near the boundary and to the well-documented instabilities of mixed-derivative schemes. As a check on the construction, the radially symmetric limit \(r_{1} \rightarrow 0\), \(\varphi_{2} - \varphi_{1} = 2\pi\) reproduces the survival probability of \citet{grebenkov2015} for the disk, the two-radius determinant collapsing to the vanishing of a single confluent hypergeometric function and the sine basis being replaced by the periodic angular basis, whose rotationally symmetric mode carries the radial result.

The paper is organized as follows. Section 2 introduces the isotropic process on the annular sector and proves the two central results: the survival probability (\S{}2.1) and the joint law of the exit time and the exit boundary, stated as a discounted boundary functional (\S{}2.2). Section 3 derives the exit-time density, its moments and Laplace transform, and resolves the joint law into boundary-piece exit densities and hitting probabilities. Section 4 establishes the companion result for a correlated, reversible OU process on a principal-axis rectangle. Section 5 describes the numerical implementation and validates the expansions against Monte Carlo simulation, including the recovery of the disk as a limiting case. Section 6 develops an application to an optically trapped colloidal particle, a setting in which the annular geometry arises naturally, and Section 7 concludes.

\section{The isotropic Ornstein--Uhlenbeck process on an annular sector}\label{sec:2}

Throughout, \(X = \left( X_{t} \right)_{t \geq 0}\) denotes the planar isotropic Ornstein--Uhlenbeck process

\begin{equation}
dX_{t} = - \beta\, X_{t}\, dt + \sigma\, dW_{t},\quad\quad X_{t} \in \mathbb{R}^{2},\ \beta > 0,\ \sigma > 0,
\tag{2.1}
\end{equation}

where \(W\) is a two-dimensional standard Brownian motion. A non-zero long-run mean is accommodated by placing the origin at that mean; as will be seen, separability requires the center of the trap and the center of the domain to coincide, so we take both at the origin without loss of generality. The generator of \(X\) is

\begin{equation}
\mathcal{L =}\frac{\sigma^{2}}{2}\,\nabla^{2} - \beta\, x \cdot \nabla,
\tag{2.2}
\end{equation}

and \(X\) is reversible with respect to the unnormalized Gaussian weight \(\rho(x) = e^{- \alpha|x|^{2}}\), \(\alpha: = \beta/\sigma^{2}\), in the sense that \(\mathcal{L =}\frac{\sigma^{2}}{2}\, e^{\alpha|x|^{2}}\nabla \cdot \left( e^{- \alpha|x|^{2}}\nabla\, \cdot \, \right)\), so that \(\mathcal{L}\) is self-adjoint on \(L^{2}(D,\rho\, dx)\) for any domain \(D\). We write the inner product as \(\langle f,g\rangle_{\rho} = \int_{D}^{}fg\,\rho\, r\, dr\, d\varphi\).

The domain is the open annular sector

\begin{equation}
D = \{(r,\varphi):\ r_{1} < r < r_{2},\ \varphi_{1} < \varphi < \varphi_{2}\},\quad\quad 0 < r_{1} < r_{2} < \infty,\ 0 \leq \varphi_{1} < \varphi_{2} \leq 2\pi,
\tag{2.3}
\end{equation}

with boundary pieces \(\Gamma_{in} = \{ r = r_{1}\}\), \(\Gamma_{out} = \{ r = r_{2}\}\) (the inner and outer arcs) and \(\Gamma_{1} = \{\varphi = \varphi_{1}\}\), \(\Gamma_{2} = \{\varphi = \varphi_{2}\}\) (the two radial edges). The boundary is \(\partial D = \Gamma_{in} \cup \Gamma_{out} \cup \Gamma_{1} \cup \Gamma_{2}\), with closure \(\overline{D} = D \cup \partial D\).

Figure 1 illustrates the annular-sector geometry and a representative sample of exit trajectories.

\begin{figure}[htbp]
\centering
\includegraphics[width=\textwidth]{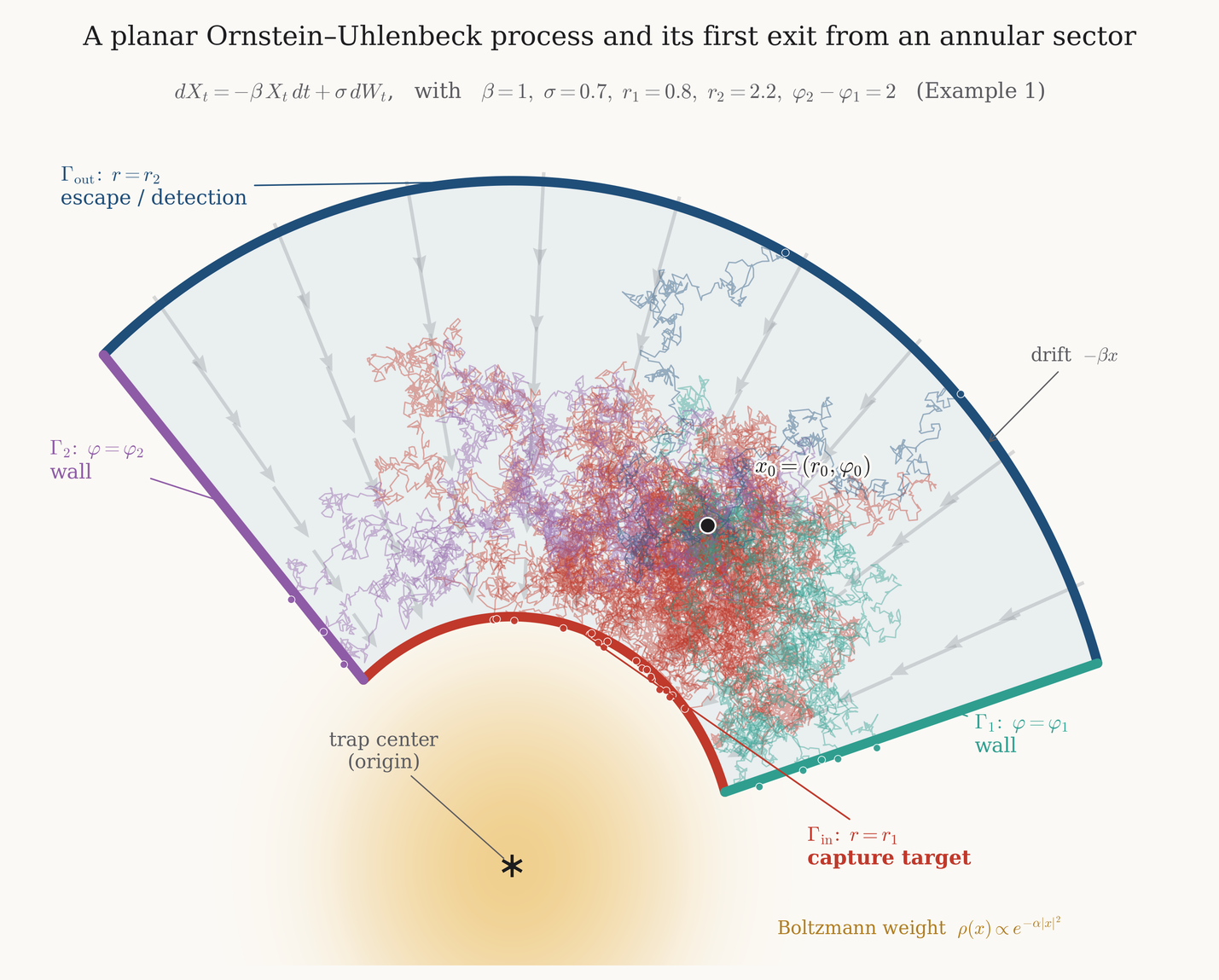}
\caption{\emph{The planar isotropic Ornstein--Uhlenbeck process and its first exit from an annular sector (the first example).} The process \(dX_{t} = - \beta\, X_{t}\, dt + \sigma\, dW_{t}\) relaxes toward the trap center at the origin, whose Gaussian equilibrium weight \(\rho(x) \propto e^{- \alpha|x|^{2}}\) is rendered as the shaded glow; arrows show the mean-reverting drift \(- \beta x\). The domain \(D\) is the open annular sector \(r_{1} < r < r_{2}\), \(\varphi_{1} < \varphi < \varphi_{2}\), with the four boundary pieces color-coded consistently throughout the paper: the inner arc \(\Gamma_{in}\) (\(r = r_{1}\)), the outer arc \(\Gamma_{out}\) (\(r = r_{2}\)), and the radial edges \(\Gamma_{1}\) (\(\varphi = \varphi_{1}\)) and \(\Gamma_{2}\) (\(\varphi = \varphi_{2}\)). Thirty representative trajectories started at \(x_{0} = \left( r_{0},\varphi_{0} \right)\) are colored by the boundary through which they exit; the visual preponderance of inner-arc exits reflects the hitting probability \(\pi_{in} \approx 0.90\). Parameters: \(\beta = 1\), \(\sigma = 0.7\), \(r_{1} = 0.8\), \(r_{2} = 2.2\), \(\varphi_{2} - \varphi_{1} = 2\), \(x_{0} = (1.30,1.00)\).}\label{fig1}
\end{figure}

Let

\begin{equation}
\tau = \inf\{ t \geq 0:\ X_{t} \notin D\}.
\tag{2.4}
\end{equation}

In polar coordinates \(x_{1} = r\cos\varphi\), \(x_{2} = r\sin\varphi\), the Laplacian is \(\nabla^{2} = \partial_{rr} + r^{- 1}\partial_{r} + r^{- 2}\partial_{\varphi\varphi}\) and the dilation field is \(x \cdot \nabla = r\partial_{r}\), so

\begin{equation}
\mathcal{L =}\frac{\sigma^{2}}{2}\left( \partial_{rr} + \frac{1}{r}\partial_{r} + \frac{1}{r^{2}}\partial_{\varphi\varphi} \right) - \beta r\,\partial_{r}.
\tag{2.5}
\end{equation}

We denote by \(M(a,b,z) =_{1}F_{1}(a;b;z)\) and \(U(a,b,z)\) the Kummer and Tricomi confluent hypergeometric functions, the two standard solutions of \(zw'' + (b - z)w' - aw = 0\) \citep[\S13.10]{olver2010}.

\subsection{Survival probability}\label{sec:2.1}

\textbf{Proposition 1.} \emph{The probability that} \(X\)\emph{, started at} \(\mathbf{x} = (r,\varphi) \in D\)\emph{, has not left} \(D\) \emph{by time} \(t\) \emph{is}

\begin{equation}
\mathbb{P}_{\mathbf{x}}(\tau > t) = \sum_{\substack{m \geq 1 \\ m\ odd}}^{}\ \sum_{p \geq 1}^{}c_{m,p}\, e^{- \mu_{m,p}\, t}\,\sin\left( \nu_{m}\left( \varphi - \varphi_{1} \right) \right)\,\mathcal{R}_{m,p}(r)
\tag{2.6}
\end{equation}

\emph{where the angular orders are the generally non-integral numbers}

\begin{equation}
\nu_{m} = \frac{m\pi}{\varphi_{2} - \varphi_{1}},\quad\quad m = 1,2,3,\ldots
\tag{2.7}
\end{equation}

\emph{the radial eigenfunctions are}

\begin{equation}
\begin{aligned}
\mathcal{R}_{m,p}(r) = \left( \alpha r^{2} \right)^{\nu_{m}/2}\Bigl[ \, & M\left( a_{m,p},\nu_{m} + 1,z_{1} \right)\, U\left( a_{m,p},\nu_{m} + 1,\alpha r^{2} \right) \\
&{}- U\left( a_{m,p},\nu_{m} + 1,z_{1} \right)\, M\left( a_{m,p},\nu_{m} + 1,\alpha r^{2} \right) \Bigr],
\end{aligned}
\tag{2.8}
\end{equation}

\emph{with} \(z_{i}: = \alpha r_{i}^{2}\) \((i = 1,2)\)\emph{, the parameter} \(a_{m,p} = \frac{1}{2}\nu_{m} - \mu_{m,p}/(2\beta)\)\emph{, and the decay rates} \(\mu_{m,p} = \beta\,\left( \nu_{m} - 2a_{m,p} \right)\) \emph{given by the increasing sequence of roots} \(a_{m,1} > a_{m,2} > \cdots \rightarrow - \infty\) \emph{of the two-radius transcendental equation}

\begin{equation}
M\left( a,\nu_{m} + 1,z_{1} \right)\, U\left( a,\nu_{m} + 1,z_{2} \right) - M\left( a,\nu_{m} + 1,z_{2} \right)\, U\left( a,\nu_{m} + 1,z_{1} \right) = 0.\
\tag{2.9}
\end{equation}

\emph{The coefficients are} \(c_{m,p} = \langle 1,\phi_{m,p}\rangle_{\rho}/ \parallel \phi_{m,p} \parallel_{\rho}^{2}\)\emph{, with} \(\phi_{m,p}(r,\varphi) = \sin\left( \nu_{m}\left( \varphi - \varphi_{1} \right) \right)\,\mathcal{R}_{m,p}(r)\)\emph{, explicitly}

\begin{equation}
c_{m,p} = \frac{2}{\nu_{m}} \cdot \frac{1}{\frac{1}{2}\left( \varphi_{2} - \varphi_{1} \right)} \cdot \frac{\int_{r_{1}}^{r_{2}}\mathcal{R}_{m,p}(r)\, e^{- \alpha r^{2}}r\, dr}{\int_{r_{1}}^{r_{2}}\mathcal{R}_{m,p}(r)^{2}\, e^{- \alpha r^{2}}r\, dr}.
\tag{2.10}
\end{equation}

\textbf{Proof.} By the connection between diffusions and parabolic equations \citep{karatzas1991,stroock1979}, \(u\left( \mathbf{x},t \right) = \mathbb{P}_{\mathbf{x}}(\tau > t)\) is the unique bounded solution of the backward initial--boundary value problem

\begin{equation}
\partial_{t}u\mathcal{= L}u\quad\text{on }D,\quad\quad u( \cdot ,0) = 1,\quad\quad u|_{\partial D} = 0\ (t > 0).
\tag{2.11}
\end{equation}

\emph{Ground-state reduction.} The first-order term in (2.5) is removed by the similarity transformation

\begin{equation}
u(r,\varphi,t) = e^{\alpha r^{2}/2}\,\psi(r,\varphi,t).
\tag{2.12}
\end{equation}

Writing \(\Phi = \alpha r^{2}\) and using \(\mathcal{L =}\frac{\sigma^{2}}{2}(\Delta - \nabla\Phi \cdot \nabla)\), a direct computation gives \(e^{- \Phi/2}\mathcal{L}\left( e^{\Phi/2}\psi \right) = \frac{\sigma^{2}}{2}\Delta\psi - V\psi\) with \(V = \frac{\sigma^{2}}{2}\left( \frac{1}{4}|\nabla\Phi|^{2} - \frac{1}{2}\Delta\Phi \right)\). Since \(|\nabla\Phi|^{2} = 4\alpha^{2}r^{2}\) and \(\Delta\Phi = 4\alpha\),

\begin{equation}
V(r) = \frac{\beta^{2}}{2\sigma^{2}}\, r^{2} - \beta,
\tag{2.13}
\end{equation}

so (2.11) becomes the harmonic-oscillator heat problem

\begin{equation}
\partial_{t}\psi = \frac{\sigma^{2}}{2}\left( \psi_{rr} + \frac{1}{r}\psi_{r} + \frac{1}{r^{2}}\psi_{\varphi\varphi} \right) - \left( \frac{\beta^{2}}{2\sigma^{2}}r^{2} - \beta \right)\psi,
\tag{2.14}
\end{equation}

\[\quad\quad\psi( \cdot ,0) = e^{- \alpha r^{2}/2},\quad\psi|_{\partial D} = 0.\quad\quad\]

\emph{Separation.} Seeking spatial eigenfunctions \(\Psi = R(r)\Theta(\varphi)\) of the operator on the right of (2.14) with eigenvalue \(- \mu\), and isolating the angular part, one obtains \(\Theta'' + \nu^{2}\Theta = 0\) on \(\left( \varphi_{1},\varphi_{2} \right)\) with \(\Theta\left( \varphi_{1} \right) = \Theta\left( \varphi_{2} \right) = 0\), whose eigenpairs are (2.7) and \(\Theta_{m}(\varphi) = \sin\left( \nu_{m}\left( \varphi - \varphi_{1} \right) \right)\).

\emph{Radial equation.} It is more transparent to separate the \emph{original} generator (2.5), \(\mathcal{L}\phi = - \mu\phi\) with \(\phi = \Theta(\varphi)R(r)\), which removes the Gaussian factor from the eigenfunctions (the eigenfunctions of \(\mathcal{L}\) on \(\mathbb{R}^{2}\) are the Hermite polynomials, not Gaussian-weighted). This yields

\begin{equation}
R'' + \left( \frac{1}{r} - 2\alpha r \right)R' + \left( \frac{2\mu}{\sigma^{2}} - \frac{\nu^{2}}{r^{2}} \right)R = 0.
\tag{2.15}
\end{equation}

The substitution \(z = \alpha r^{2}\), \(R = z^{\nu/2}v(z)\) turns (2.15) into Kummer's equation

\begin{equation}
z\, v'' + (\nu + 1 - z)\, v' - a\, v = 0,\quad\quad a = \frac{\nu}{2} - \frac{\mu}{2\beta},
\tag{2.16}
\end{equation}

equivalently \(\mu = \beta(\nu - 2a)\). (As a check, the polynomial --- free-space --- solutions \(a = - n_{r}\), \(n_{r} \in \mathbb{Z}_{\geq 0}\), give \(\mu = \beta\left( \nu + 2n_{r} \right)\); for integer \(\nu = \left| \ell \right|\) this is the known OU spectrum \(\mu \in \beta\,\mathbb{Z}_{\geq 0}\).) The general radial solution is therefore

\begin{equation}
R(r) = \left( \alpha r^{2} \right)^{\nu/2}\left\lbrack c_{1}M\left( a,\nu + 1,\alpha r^{2} \right) + c_{2}U\left( a,\nu + 1,\alpha r^{2} \right) \right\rbrack.
\tag{2.17}
\end{equation}

\emph{Eigenvalues.} The Dirichlet conditions \(R\left( r_{1} \right) = R\left( r_{2} \right) = 0\) give a homogeneous linear system in \(\left( c_{1},c_{2} \right)\); the common prefactor \(\left( \alpha r_{i}^{2} \right)^{\nu/2}\) cancels, and a non-trivial solution exists iff the determinant vanishes, which is (2.9). The determinant (2.9) vanishes spuriously at every non-positive integer \(a \in \{ 0, - 1, - 2,\ldots\}\), where the Kummer and Tricomi functions become proportional --- both reducing to the same generalized Laguerre polynomial --- so that (2.9) vanishes there for want of two independent radial solutions rather than because an eigenfunction exists. Generically this polynomial does not vanish at both radii, no radial eigenfunction meets the two boundary conditions, and these roots are discarded; only in the non-generic case that it vanishes simultaneously at both radii is such a value a genuine eigenvalue, recovered as the limit of the neighboring non-integer roots. The eigenvalues \(\mu_{m,p}\) are obtained from the non-integer roots. Equation (2.15) is, after multiplication by the integrating factor \(w(r) = r\, e^{- \alpha r^{2}}\), in the regular Sturm--Liouville form \((wR')' - \frac{\nu^{2}}{r^{2}}wR + \frac{2\mu}{\sigma^{2}}wR = 0\) on \(\left\lbrack r_{1},r_{2} \right\rbrack\); because \(r_{1} > 0\) the weight \(w\) is bounded and bounded away from \(0\) and both endpoints are \emph{regular}, so for each \(m\) the spectrum \(\{\mu_{m,p}\}_{p \geq 1}\) is real, simple, and increases to \(+ \infty\) \citep{zettl2005}, and the family \(\{\phi_{m,p}\}_{m,p}\) is complete and orthogonal in \(L^{2}(D,\rho)\). Taking \(c_{1} = - U\left( a,\nu + 1,z_{1} \right)\), \(c_{2} = M\left( a,\nu + 1,z_{1} \right)\) produces the eigenfunction (2.8), which vanishes at \(r_{1}\) by construction and at \(r_{2}\) by (2.9).

\emph{Assembly.} Expanding the initial datum \(1 \in L^{2}(D,\rho)\) in this basis,

\begin{equation}
u\left( \mathbf{x},t \right) = \sum_{m,p}^{}c_{m,p}\, e^{- \mu_{m,p}t}\,\phi_{m,p}\left( \mathbf{x} \right),\quad\quad c_{m,p} = \frac{\langle 1,\phi_{m,p}\rangle_{\rho}}{\parallel \phi_{m,p} \parallel_{\rho}^{2}},
\tag{2.18}
\end{equation}

convergent in \(L^{2}(\rho)\) and, by parabolic regularity, classically for \(t > 0\). The angular factor of the numerator is \(\int_{\varphi_{1}}^{\varphi_{2}}\sin\left( \nu_{m}\left( \varphi - \varphi_{1} \right) \right)\, d\varphi = \nu_{m}^{- 1}\left( 1 - ( - 1)^{m} \right)\), which vanishes for even \(m\) and equals \(2/\nu_{m}\) for odd \(m\); the denominator's angular factor is \(\int_{\varphi_{1}}^{\varphi_{2}}\sin^{2}\left( \nu_{m}\left( \varphi - \varphi_{1} \right) \right)\, d\varphi = \frac{1}{2}\left( \varphi_{2} - \varphi_{1} \right)\). This gives (2.10) and the restriction to odd \(m\) in (2.6). \(\quad\quad \ensuremath{\square}\)

\textbf{Remark 1 (radial quadratures).} The two radial integrals in (2.10) are, under \(z = \alpha r^{2}\), of the form \(\int_{z_{1}}^{z_{2}}z^{c}e^{- z}\, M(a,b,z)\, dz\) and \(\int_{z_{1}}^{z_{2}}z^{c}e^{- z}\, U(a,b,z)\, dz\), which are evaluated either by stable one-dimensional quadrature or through antiderivative identities for confluent hypergeometric functions of shifted parameters \citep[\S13.10]{olver2010}; alternatively, the squared norm \(\int_{r_{1}}^{r_{2}}\mathcal{R}_{m,p}^{2}\, e^{- \alpha r^{2}}r\, dr\) is obtained in closed form from the Sturm--Liouville normalization as a boundary term, namely \(- \frac{\sigma^{2}}{2}\, w(r)\,\mathcal{R}_{m,p}(r)\,\partial_{\mu}\mathcal{R}_{m,p}'(r)|_{r_{1}}^{r_{2}}\) evaluated along (2.9).

\textbf{Remark 2 (Whittaker form).} Writing the radial solutions through Whittaker functions, \(\mathcal{R}_{m,p}(r) = \frac{e^{\alpha r^{2}/2}}{\sqrt{\alpha}r}\left\lbrack c_{1}M_{\kappa_{m,p},\nu_{m}/2}\left( \alpha r^{2} \right) + c_{2}W_{\kappa_{m,p},\nu_{m}/2}\left( \alpha r^{2} \right) \right\rbrack\) with \(\kappa_{m,p} = \frac{1}{2} + \mu_{m,p}/(2\beta)\) and second index \(\nu_{m}/2\), exhibits (2.8) as a two-radius combination of Whittaker functions of generally non-integral order.

\textbf{Remark 3 (the disk limit and Grebenkov).} Letting \(r_{1} \rightarrow 0\) and \(\varphi_{2} - \varphi_{1} \rightarrow 2\pi\) degenerates the regular Sturm--Liouville problem into a \emph{singular} one: the inner endpoint becomes a singular point of (2.15), the irregular solution \(U\) (which behaves like \(z^{- \nu/2}\), i.e.~\(r^{- \nu}\)) is excluded by the limit-point condition at the origin, and (2.9) collapses to the single-function condition \(M\left( a,1,z_{2} \right) = 0\) on the rotationally symmetric mode \(\nu = 0\). One then recovers the radially symmetric survival probability of a harmonically trapped particle in a disk obtained by Grebenkov \citeyearpar{grebenkov2015} as zeros of a confluent hypergeometric function. This collapse of the radial problem requires the angular boundary conditions to be changed as well: under the absorbing radial edges retained here the angular orders are the generally non-integral numbers \(\frac{m\pi}{\varphi_{2} - \varphi_{1}}\), so that as \(\varphi_{2} - \varphi_{1} \rightarrow 2\pi\) the lowest of them tends to \(\frac{1}{2}\), never to the rotationally symmetric value \(0\). The disk and full annulus are therefore not literal Dirichlet-sector limits; they are recovered from the same radial Kummer equation once the sine basis is replaced by the periodic angular basis, whose zero mode reproduces the radial result. The strict positivity of \(r_{1}\) is precisely what restores a regular eigenproblem requiring both confluent hypergeometric solutions, and herein lies the structural novelty of the annular geometry: in exactly-soluble whole-plane ring models such as that of Tan and Inkson \citeyearpar{tan1996}, the origin belongs to the domain, so --- as above --- the irregular Tricomi solution is excluded there and a single Kummer function fixes the spectrum, whereas the inner boundary here retains both solutions and yields the two-radius determinant (2.9).

\subsection{Joint law of the exit time and the exit boundary}\label{sec:2.2}

The object of primary interest is the joint distribution of \(\left( \tau,X_{\tau} \right)\) --- in particular, \emph{which} of the four boundary pieces is reached and after how long. For a radially symmetric domain this refinement is empty; here it engages the full eigenbasis of Proposition 1, including the even angular modes absent from (2.6). We state it in the form of a discounted boundary functional from which the time-resolved joint law and the boundary-hitting probabilities follow as special cases.

\textbf{Proposition 2.} \emph{Let} \(\varrho \geq 0\) \emph{be a discount rate and let} \(g_{in},g_{out}:\lbrack 0\mathbb{,\infty) \rightarrow R}\) \emph{and} \(h_{1},h_{2}:\lbrack 0\mathbb{,\infty) \rightarrow R}\) \emph{be piecewise-continuous functions, interpreted as the amount paid if} \(X\) \emph{first leaves} \(D\) \emph{through} \(\Gamma_{in}\)\emph{,} \(\Gamma_{out}\)\emph{,} \(\Gamma_{1}\)\emph{,} \(\Gamma_{2}\) \emph{respectively. Then}

\begin{equation}
\begin{aligned}
\Pi\left( \mathbf{x},t \right): = \mathbb{E}_{\mathbf{x}}\biggl[ e^{- \varrho\tau}\Bigl( & g_{in}(\tau)\mathbf{1}_{X_{\tau} \in \Gamma_{in}} + g_{out}(\tau)\mathbf{1}_{X_{\tau} \in \Gamma_{out}} \\
&{}+ h_{1}(\tau)\mathbf{1}_{X_{\tau} \in \Gamma_{1}} + h_{2}(\tau)\mathbf{1}_{X_{\tau} \in \Gamma_{2}} \Bigr)\,\mathbf{1}_{\tau \leq t} \biggr]
\end{aligned}
\tag{2.19}
\end{equation}

\emph{carries, when the four rewards are independent of time, the boundary values} \(\Pi = g_{in}\) \emph{on} \(\Gamma_{in}\)\emph{,} \(\Pi = g_{out}\) \emph{on} \(\Gamma_{out}\)\emph{,} \(\Pi = h_{1}\) \emph{on} \(\Gamma_{1}\)\emph{,} \(\Pi = h_{2}\) \emph{on} \(\Gamma_{2}\)\emph{, and is in every case given by the boundary-flux representation}

\begin{equation}
\Pi\left( \mathbf{x},t \right) = \frac{\sigma^{2}}{2}\int_{0}^{t}e^{- \varrho s}\int_{\partial D}^{}\rho(\xi)\,\partial_{n_{\xi}}G\left( \mathbf{x},\xi;s \right)\, H_{\partial D}(\xi,s)\, d\ell(\xi)\, ds,
\tag{2.20}
\end{equation}

\emph{where} \(n_{\xi}\) \emph{is the inward normal,} \(d\ell\) \emph{is arclength,} \(\rho\) \emph{is the Gaussian weight introduced in Section 2,} \(H_{\partial D}\) \emph{denotes the boundary data above on the piece containing} \(\xi\)\emph{, evaluated at the exit time, and} \(G\) \emph{is the Dirichlet heat kernel of} \(\mathcal{L}\) \emph{on} \(D\)\emph{,}

\begin{equation}
G\left( \mathbf{x},\xi;s \right) = \sum_{m \geq 1}^{}{\sum_{p \geq 1}^{}e^{- \mu_{m,p}s}}\,\frac{\phi_{m,p}\left( \mathbf{x} \right)\,\phi_{m,p}(\xi)}{\parallel \phi_{m,p} \parallel_{\rho}^{2}}
\tag{2.21}
\end{equation}

\emph{the eigenpairs being those of Proposition 1.}

\emph{When the four payoffs are independent of time, the boundary reward is a function of position alone and} \(\Pi\) \emph{is the unique bounded solution of the homogeneous parabolic equation} \(\partial_{t}\Pi = \left( \mathcal{L -}\varrho \right)\Pi\) \emph{on} \(D\)\emph{, with zero initial data} \(\Pi( \cdot ,0) = 0\) \emph{and boundary data} \(\Pi = H_{\partial D}\) \emph{on} \(\partial D\)\emph{; (2.20) is then its Duhamel representation, and its infinite-horizon limit solves the corresponding elliptic boundary-value problem with the same data. When a payoff depends on the exit time, the corresponding boundary value inherits that dependence, an inhomogeneous boundary value on} \(\partial D\) \emph{enters the parabolic Poisson representation, and} \(\Pi\) \emph{is characterized by the boundary representation (2.20) rather than by a homogeneous equation in} \(t\)\emph{.}

\emph{In particular:} (i) \emph{with} \(\varrho = 0\) \emph{and} \(g_{out} \equiv \mathbf{1}\)\emph{, all other data zero,} \(\Pi\left( \mathbf{x},t \right) = \mathbb{P}_{\mathbf{x}}\left( \tau \leq t,\, X_{\tau} \in \Gamma_{out} \right)\) \emph{is the probability of leaving through the outer arc by time} \(t\)\emph{, and analogously for each piece;} (ii) \emph{the four integrands} \(\frac{\sigma^{2}}{2}\rho(\xi)\,\,\partial_{n}G\left( \mathbf{x},\xi;s \right)\)\emph{, restricted to} \(\xi \in \Gamma_{in},\Gamma_{out},\Gamma_{1},\Gamma_{2}\)\emph{, are the joint sub-densities of} \(\left( \tau,X_{\tau} \right)\) \emph{on the respective pieces; integrating the outer-arc and inner-arc densities over their arcs and the edge densities over their edges, and summing, returns} \(- \partial_{t}\mathbb{P}_{\mathbf{x}}(\tau > t)\) \emph{from (2.6). Allowing the payoffs to depend on the exit position as well as the exit time changes nothing in the derivation: the boundary data in (2.20) are then read as piecewise-continuous functions on the boundary itself, and the integrands in (ii) resolve the joint law jointly in the exit time and the exit position.}

\textbf{Proof.} The discount contributes only the multiplicative factor \(\exp( - \varrho\tau)\), equal to \(\exp( - \varrho s)\) at the exit instant \(\tau = s\); it is therefore enough to establish the representation in the undiscounted case \(\varrho = 0\) and to restore this factor under the time integral, which is the effect of the change of unknown that removes the zeroth-order term. By the correspondence between the diffusion and the Dirichlet heat equation already used for Proposition 1, the sub-Markovian semigroup killed on \(\partial D\) acts, for each \(s > 0\), as an integral operator against \(\rho\, d\xi\),

\[E_{x}\left\lbrack f\left( X_{s} \right)\, 1_{s < \tau} \right\rbrack = \int_{D}^{}{G(x,\xi;s)\, f(\xi)\,\rho(\xi)\, d\xi}\]

its kernel \(G\) being the Dirichlet heat kernel of \(L\) on \(D\). The operator \(- L\) is self-adjoint and has compact resolvent on \(L^{2}(D,\rho\, dx)\) under the absorbing conditions of Proposition 1, with eigenpairs \(\left( \mu_{m,p},\phi_{m,p} \right)\) and squared norms \(\left\| \phi_{m,p} \right\|^{2}\); the spectral theorem yields the bilinear expansion (2.21), convergent for every \(s > 0\) and absolutely once \(s\) is bounded away from \(0\). Equivalently, \(G( \cdot ,\xi;s)\) is the unique bounded solution of \(\partial_{s}G = LG\) in \(D\) that vanishes on \(\partial D\) and tends to \(\delta_{\xi}\), relative to \(\rho\, d\xi\), as \(s \downarrow 0\); the kernel is symmetric, \(G(x,\xi;s) = G(\xi,x;s)\), and solves the heat equation in either argument.

The reversibility of the process puts \(L\) in divergence form, \(L = \left( \frac{\sigma^{2}}{2} \right)\rho^{- 1}\nabla \cdot (\rho\nabla \cdot )\), so that Green's identity carries the weight \(\rho\):

\[\int_{D}^{}{(Lu)v\rho\, dx} = - \frac{\sigma^{2}}{2}\int_{D}^{}{\rho\nabla u \cdot \nabla v\, dx} + \frac{\sigma^{2}}{2}\oint_{\partial D}^{}{\rho v\partial_{\nu}u\, dl}\]

with \(\partial_{\nu}\) the outward normal derivative and \(\, dl\) arclength. Taking \(u = G(x, \cdot ;s)\) and \(v \equiv 1\), and using \(\partial_{s}G = L_{\xi}G\), one differentiates the survival probability \(P_{x}(\tau > s) = \int_{D}^{}{G(x,\xi;s)\rho(\xi)\, d\xi}\) of (2.6) to obtain

\[- \frac{d}{ds}\, P_{x}(\tau > s) = \frac{\sigma^{2}}{2}\oint_{\partial D}^{}{\rho(\xi)\,\partial_{n(\xi)}\, G(x,\xi;s)\, dl(\xi)}\]

where \(n = - \nu\) is the inward normal; the integrand is non-negative, since \(G\) is positive in \(D\) and vanishes on \(\partial D\), and the left-hand side is the total rate of exit at time \(s\). Replacing the constant \(1\) by a function supported near a single boundary piece localizes the identity and shows that the pair \(\left( \tau,X_{\tau} \right)\) possesses, on \((0,\infty) \times \partial D\), the joint density \(p(x;s,\xi) = \left( \frac{\sigma^{2}}{2} \right)\rho(\xi)\,\partial_{n}(\xi)G(x,\xi;s)\) with respect to \(\, ds\) and arclength \(\, dl(\xi)\). The weight \(\rho\) is exactly the one under which \(L\) is self-adjoint, and it is the factor carried by the key density (3.7).

The eigenfunctions of Proposition 1 are \(\phi_{m,p}(r,\varphi) = \sin\left( \nu_{m}\left( \varphi - \varphi_{1} \right) \right)R_{m,p}(r)\), so the inward normal derivative of \(G\) is obtained termwise from (2.21): on the arcs \(\Gamma_{in}\) and \(\Gamma_{out}\) it equals \(\mp \partial_{r}\phi_{m,p}\), while on the radial edges \(\Gamma_{1}\) and \(\Gamma_{2}\) it equals \(\mp r^{- 1}\partial_{\varphi}\phi_{m,p} = \mp \nu_{m}\, r^{- 1}\cos\left( \nu_{m}\left( \varphi - \varphi_{1} \right) \right)R_{m,p}(r)\), which does not vanish at \(\varphi = \varphi_{1}\) and \(\varphi = \varphi_{2}\) and therefore brings in the even angular modes that the symmetric datum of Proposition 1 annihilated; these are the series that enter Corollary 3.

Since the four events \(\left\{ X_{\tau} \in \Gamma_{in} \right\}\), \(\left\{ X_{\tau} \in \Gamma_{out} \right\}\), \(\left\{ X_{\tau} \in \Gamma_{1} \right\}\), \(\left\{ X_{\tau} \in \Gamma_{2} \right\}\) partition \(\left\{ \tau < \infty \right\}\) up to a null set, and the boundary data are piecewise continuous and bounded on compact time intervals, the joint density together with Fubini's theorem gives

\[\Pi(x,t) = E_{x}\left\lbrack e^{- \varrho\tau}\, H_{\partial D}\left( X_{\tau},\tau \right)\, 1_{\tau \leq t} \right\rbrack = \int_{0}^{t}{\oint_{\partial D}^{}{e^{- \varrho s}\, H_{\partial D}(\xi,s)\, p(x;s,\xi)\, dl(\xi)\, ds}}\]

which is (2.20) after substitution of the joint density; here \(H_{\partial D}(\xi,s)\) denotes \(g_{in}(s)\), \(g_{out}(s)\), \(h_{1}(s)\), \(h_{2}(s)\) on the respective pieces, evaluated at the exit time \(s\). When the data do not depend on time, the same identity characterizes the stationary (infinite-horizon) limit of \(\Pi\) as the bounded solution of the elliptic problem \((L - \varrho)\Pi = 0\) in \(D\) with \(\Pi = H_{\partial D}\) on \(\partial D\), the maximum-principle characterization whose uniqueness is classical. Taking \(H_{\partial D}\) to be the indicator of a single piece, with \(\varrho = 0\), gives the boundary-hitting probability of statement (i); leaving the data general and reading off the integrand gives the joint sub-densities of statement (ii), whose mutual consistency --- that they integrate over \(\partial D\) to \(- \frac{d}{ds}P_{x}(\tau > s)\) --- is precisely the flux identity established above. \(\ \ \ensuremath{\square}\)

\section{Exit-time density, moments and the joint exit law}\label{sec:3}

Throughout this section the eigenpairs \(\left( \mu_{m,p},\phi_{m,p} \right)\), the angular orders \(\nu_{m}\), the radial eigenfunctions \(\mathcal{R}_{m,p}\) and the coefficients \(c_{m,p}\) are those of Proposition 1; we write \(\Delta\varphi: = \varphi_{2} - \varphi_{1}\), \(\rho(x) = e^{- \alpha|x|^{2}}\), and

\begin{equation}
N_{m,p}: = \parallel \phi_{m,p} \parallel_{\rho}^{2} = \frac{\Delta\varphi}{2}\,\eta_{m,p},\quad\quad\eta_{m,p}: = \int_{r_{1}}^{r_{2}}\mathcal{R}_{m,p}(r)^{2}\, e^{- \alpha r^{2}}r\, dr,
\tag{3.1}
\end{equation}

together with the two radial functionals that will carry all the explicit content,

\begin{equation}
\mathcal{K}_{m,p}: = \int_{r_{1}}^{r_{2}}\mathcal{R}_{m,p}(r)\, e^{- \alpha r^{2}}r\, dr,\quad\quad\mathcal{J}_{m,p}: = \int_{r_{1}}^{r_{2}}\mathcal{R}_{m,p}(r)\,\frac{e^{- \alpha r^{2}}}{r}\, dr.
\tag{3.2}
\end{equation}

Both are linear confluent-hypergeometric moment integrals with closed primitives \citep[\S13.10]{olver2010} (the squared norms, by contrast, are obtained by quadrature or the Sturm--Liouville normalization of Remark 1), so that, for odd \(m\), \(c_{m,p} = 4\mathcal{K}_{m,p}/\left( \nu_{m}\,\Delta\varphi\,\eta_{m,p} \right)\).

\subsection{Density and moments of the exit time}\label{sec:3.1}

\textbf{Corollary 1.} \emph{The first-exit time} \(\tau\) \emph{from} \(D\) \emph{has, for} \(\mathbf{x} = (r,\varphi) \in D\)\emph{, the density}

\begin{equation}
f_{\tau}\left( \mathbf{x},t \right) = - \,\partial_{t}\,\mathbb{P}_{\mathbf{x}}(\tau > t) = \sum_{\substack{m \geq 1 \\ m\ odd}}^{}\ \sum_{p \geq 1}^{}c_{m,p}\,\mu_{m,p}\, e^{- \mu_{m,p}t}\,\sin\left( \nu_{m}\left( \varphi - \varphi_{1} \right) \right)\,\mathcal{R}_{m,p}(r).
\tag{3.3}
\end{equation}

\emph{As} \(t \rightarrow \infty\) \emph{the density is governed by the smallest decay rate occurring in (3.3), namely} \(\mu_{*}: = \mu_{1,1}\) \emph{(the lowest radial mode of the first angular harmonic), and}

\begin{equation}
f_{\tau}\left( \mathbf{x},t \right)\  \sim \ c_{1,1}\,\mu_{1,1}\, e^{- \mu_{1,1}t}\,\sin\left( \nu_{1}\left( \varphi - \varphi_{1} \right) \right)\,\mathcal{R}_{1,1}(r),\quad\quad t \rightarrow \infty.
\tag{3.4}
\end{equation}

The proof is immediate from termwise differentiation of (2.6), legitimate for \(t > 0\) by parabolic regularity, the differentiated series converging locally uniformly on compact subsets of \(D\). We record the spectral expressions for the moments and the Laplace transform.

\textbf{Corollary 2.} \emph{For every} \(\mathbf{x} \in D\) \emph{and every integer} \(n \geq 1\)\emph{,}

\begin{equation}
\mathbb{E}_{\mathbf{x}}\left\lbrack \tau^{\, n} \right\rbrack = n!\sum_{m\ odd}^{}{\sum_{p \geq 1}^{}\frac{c_{m,p}}{\mu_{m,p}^{\, n}}}\,\sin\left( \nu_{m}\left( \varphi - \varphi_{1} \right) \right)\,\mathcal{R}_{m,p}(r),
\tag{3.5}
\end{equation}

\emph{in particular the mean exit time} \(w\left( \mathbf{x} \right) = \mathbb{E}_{\mathbf{x}}\lbrack\tau\rbrack\) \emph{is the unique solution of the Andronov--Vitt--Pontryagin problem} \(\mathcal{L}w = - 1\) \emph{on} \(D\)\emph{,} \(w|_{\partial D} = 0\)\emph{, given by (3.5) with} \(n = 1\)\emph{. The Laplace transform of} \(\tau\) \emph{is}

\begin{equation}
\mathbb{E}_{\mathbf{x}}\left\lbrack e^{- \varrho\tau} \right\rbrack = \sum_{m\ odd}^{}{\sum_{p \geq 1}^{}\frac{\mu_{m,p}}{\varrho + \mu_{m,p}}}\, c_{m,p}\,\sin\left( \nu_{m}\left( \varphi - \varphi_{1} \right) \right)\,\mathcal{R}_{m,p}(r),\quad\quad\varrho \geq 0,
\tag{3.6}
\end{equation}

\emph{and solves} \(\left( \mathcal{L -}\varrho \right)v_{\varrho} = 0\,\) \emph{on} \(D\) \emph{with} \(v_{\varrho}|_{\partial D} = 1\)\emph{.}

\textbf{Proof.} Identities (3.5) and (3.6) follow from (3.3) and the elementary integrals \(\int_{0}^{\infty}t^{n - 1}e^{- \mu t}dt = (n - 1)!\,\mu^{- n}\) and \(\int_{0}^{\infty}e^{- \varrho t}\mu e^{- \mu t}dt = \mu/(\varrho + \mu)\), applied termwise. The probabilistic characterizations \(\mathcal{L}w = - 1\) and \(\left( \mathcal{L -}\varrho \right)v_{\varrho} = 0\) are the standard Dynkin/Feynman--Kac equations for the exit functionals; that the series solve them in the interior is verified termwise, the constant boundary datum being recovered through the \(L^{2}(\rho)\)-completeness of \(\{\phi_{m,p}\}\). Setting \(\varrho = 0\) in (3.6) returns \(\sum_{m\ odd,p}^{}c_{m,p}\phi_{m,p} = 1\), i.e.~\(\mathbb{P}_{\mathbf{x}}(\tau < \infty) = 1\). \(\quad\quad \ensuremath{\square}\)

Equation (3.5) with \(n = 1\) extends to the annular sector the mean-exit-time computation carried out for the disk and the ball by Grebenkov \citeyearpar{grebenkov2015} and, asymptotically in high dimension, by Kersting et al.~\citeyearpar{kersting2023}; the present formula resolves the angular as well as the radial dependence.

\subsection{The joint law of the exit time and the exit boundary}\label{sec:3.2}

We now turn to the principal object: the joint distribution of \(\left( \tau,X_{\tau} \right)\), resolving \emph{which} of the four boundary pieces is reached and \emph{when}. As anticipated in Section 2, this engages the entire eigenbasis \(\{\phi_{m,p}\}_{m \geq 1,p \geq 1}\) --- including the even angular modes annihilated by the symmetric initial datum of Proposition 1 --- and constitutes the genuinely new content that the radial symmetry of the disk problem renders inaccessible.

We first record the joint exit density, the specialization of Proposition 2 to a Dirac datum in time and an indicator in space. Let \(\partial_{n}\) denote the \emph{inward} normal derivative on \(\partial D\).

\textbf{Corollary 3 (joint exit density).} \emph{The pair} \(\left( \tau,X_{\tau} \right)\) \emph{has, for} \(\mathbf{x} \in D\)\emph{, the joint density with respect to} \(dt\, d\ell(\xi)\) \emph{on} \((0,\infty) \times \partial D\)\emph{, away from the four corners, given by}

\begin{equation}
\mathfrak{j}\left( t;\mathbf{x},\xi \right) = \frac{\sigma^{2}}{2}\, e^{- \alpha r_{\xi}^{2}}\sum_{m \geq 1}^{}{\sum_{p \geq 1}^{}e^{- \mu_{m,p}t}}\,\frac{\phi_{m,p}\left( \mathbf{x} \right)\,\partial_{n}\phi_{m,p}(\xi)}{N_{m,p}},\quad\quad\xi \in \partial D
\tag{3.7}
\end{equation}

Exit at one of the four corners, where two boundary pieces meet and the outward normal is undefined, occurs with probability zero, so (3.7) is the joint density on the open arcs and open radial edges; the boundary normal derivative is taken on these smooth pieces, and the marginal sub-densities of Corollary 4 are obtained by integrating (3.7) over each of them. There, for \(t > 0\), the series is the termwise inward-normal derivative of the heat-kernel expansion (2.21), which converges once integrated over an arc or an edge, the pointwise convergence up to the boundary being slower and boundary-layer sensitive.

\textbf{Proof.} Let \(q\left( t;\mathbf{x},\mathbf{y} \right)\) be the sub-density, with respect to Lebesgue measure, of the process killed on leaving \(D\), so that \(\mathbb{P}_{\mathbf{x}}\left( X_{t} \in d\mathbf{y},\tau > t \right) = q\, d\mathbf{y}\). Writing \(\mathcal{L}\) in the divergence form \(\mathcal{L =}\frac{\sigma^{2}}{2}\,\rho^{- 1}\nabla \cdot (\rho\nabla\, \cdot \,)\) established in Section 2, self-adjointness on \(L^{2}(D,\rho)\) gives the bilinear expansion \(q\left( t;\mathbf{x},\mathbf{y} \right) = \rho\left( \mathbf{y} \right)\sum_{m,p}^{}e^{- \mu_{m,p}t}\phi_{m,p}\left( \mathbf{x} \right)\phi_{m,p}\left( \mathbf{y} \right)/N_{m,p}\), that is, \(\rho\) times the symmetric heat kernel \(G\) of (2.21) --- the symmetric kernel being taken with respect to the Gaussian weight, not to Lebesgue measure --- and its integral against \(d\mathbf{y}\) of the datum \(1\) reproduces (2.6). For a diffusion with diffusion matrix \(\sigma^{2}I\), the joint exit density is the inward probability flux \(\mathfrak{j}\left( t;\mathbf{x},\xi \right) = \frac{\sigma^{2}}{2}\,\partial_{n}q\left( t;\mathbf{x},\xi \right)\) \citep{karatzas1991,stroock1979}. Since every \(\phi_{m,p}\) vanishes on \(\partial D\), the product rule leaves only the term in which \(\partial_{n}\) falls on \(\phi_{m,p}(\xi)\), yielding (3.7). Equivalently, if \(K^{D}(t;x,y)\) denotes the symmetric Dirichlet heat kernel with respect to the weighted measure \(\rho(y)\, dy\), then the killed transition density with respect to Lebesgue measure is \(q(t;x,y) = \rho(y)K^{D}(t;x,y).\) Since \(K^{D}\), and hence \(q\), vanishes on \(\partial D\), the drift term contributes no boundary flux there, and the joint exit density with respect to time and arclength is \(j(t;x,y) = \frac{\sigma^{2}}{2}\rho(y)\,\partial_{n_{in},y}K^{D}(t;x,y),y \in \partial D.\)

Substituting the eigenfunction expansion of \(K^{D}\) gives (3.7). For verification, integrating (3.7) over \(t \in (0,\infty)\) and \(\xi \in \partial D\) and using the divergence theorem together with \(\mathcal{L}\phi_{m,p} = - \mu_{m,p}\phi_{m,p}\) gives \(\oint_{\partial D}^{}e^{- \alpha r_{\xi}^{2}}\partial_{n}\phi_{m,p}\, d\ell =\frac{2\mu_{m,p}}{\sigma^{2}}\langle 1,\phi_{m,p}\rangle_{\rho}\), whence the total mass equals \(\sum_{m,p}^{}c_{m,p}\phi_{m,p}\left( \mathbf{x} \right) = 1\). \(\quad\quad \ensuremath{\square}\)

The inward normal derivatives of \(\phi_{m,p} = \sin\left( \nu_{m}\left( \varphi - \varphi_{1} \right) \right)\mathcal{R}_{m,p}(r)\) on the four boundary pieces are elementary:

\begin{equation}
\begin{gathered}
\partial_{n}\phi_{m,p}|_{\Gamma_{out}} = - \mathcal{R}_{m,p}^{'}\left( r_{2} \right)\sin\left( \nu_{m}\left( \psi - \varphi_{1} \right) \right), \\
\partial_{n}\phi_{m,p}|_{\Gamma_{in}} = + \mathcal{R}_{m,p}'\left( r_{1} \right)\sin\left( \nu_{m}\left( \psi - \varphi_{1} \right) \right)
\end{gathered}
\tag{3.8}
\end{equation}

\begin{equation}
\partial_{n}\phi_{m,p}|_{\Gamma_{1}} = \frac{\nu_{m}}{s}\,\mathcal{R}_{m,p}(s),\quad\quad\partial_{n}\phi_{m,p}|_{\Gamma_{2}} = \frac{\nu_{m}( - 1)^{m + 1}}{s}\,\mathcal{R}_{m,p}(s)
\tag{3.9}
\end{equation}

where \(\psi \in \left( \varphi_{1},\varphi_{2} \right)\) parameterizes the arcs and \(s \in \left( r_{1},r_{2} \right)\) the radial edges. The contrast between (3.8) and (3.9) is the structural heart of the result: on the \emph{arcs} the angular factor is the sine, which integrates to zero over a full angular span for even \(m\); on the \emph{edges} the angular factor is the cosine \(\partial_{\varphi}\sin\), which does not vanish, so that \emph{all} angular modes contribute. Integrating (3.7) against the relevant boundary coordinate gives the four marginal-in-time exit densities.

\textbf{Corollary 4 (boundary-resolved exit laws).} \emph{For} \(\mathbf{x} = (r,\varphi) \in D\)\emph{, the sub-densities of} \(\tau\) \emph{on the four boundary pieces are}

\begin{equation}
f_{\tau}^{out}\left( \mathbf{x},t \right) = 2\sigma^{2}r_{2}\, e^{- \alpha r_{2}^{2}}\sum_{m\ odd}^{}{\sum_{p}^{}\frac{- \,\mathcal{R}_{m,p}'\left( r_{2} \right)}{\nu_{m}\,\Delta\varphi\,\eta_{m,p}}}\mspace{6mu} e^{- \mu_{m,p}t}\,\phi_{m,p}\left( \mathbf{x} \right)
\tag{3.10}
\end{equation}

\begin{equation}
f_{\tau}^{in}\left( \mathbf{x},t \right) = 2\sigma^{2}r_{1}\, e^{- \alpha r_{1}^{2}}\sum_{m\ odd}^{}{\sum_{p}^{}\frac{+ \,\mathcal{R}_{m,p}'\left( r_{1} \right)}{\nu_{m}\,\Delta\varphi\,\eta_{m,p}}}\mspace{6mu} e^{- \mu_{m,p}t}\,\phi_{m,p}\left( \mathbf{x} \right)
\tag{3.11}
\end{equation}

\begin{equation}
\begin{gathered}
f_{\tau}^{\, 1}\left( \mathbf{x},t \right) = \frac{\sigma^{2}}{\Delta\varphi}\sum_{m \geq 1}^{}{\sum_{p}^{}\frac{\nu_{m}\,\mathcal{J}_{m,p}}{\eta_{m,p}}}\mspace{6mu} e^{- \mu_{m,p}t}\,\phi_{m,p}\left( \mathbf{x} \right),\quad\quad \\
f_{\tau}^{\, 2}\left( \mathbf{x},t \right) = \frac{\sigma^{2}}{\Delta\varphi}\sum_{m \geq 1}^{}{\sum_{p}^{}\frac{( - 1)^{m + 1}\nu_{m}\,\mathcal{J}_{m,p}}{\eta_{m,p}}}\mspace{6mu} e^{- \mu_{m,p}t}\,\phi_{m,p}\left( \mathbf{x} \right)
\end{gathered}
\tag{3.12}
\end{equation}

\emph{The corresponding boundary-hitting probabilities} \(\pi_{\bullet}\left( \mathbf{x} \right) = \mathbb{P}_{\mathbf{x}}\left( X_{\tau} \in \Gamma_{\bullet} \right)\) \emph{are obtained by replacing each} \(e^{- \mu_{m,p}t}\) \emph{by} \(\mu_{m,p}^{- 1}\)\emph{, understood as the zero-discount limit of the discounted boundary-hitting probabilities --- a termwise time integration that is legitimate in this Abel sense even where the sub-density series converges only conditionally --- and satisfy}

\begin{equation}
\pi_{out}\left( \mathbf{x} \right) + \pi_{in}\left( \mathbf{x} \right) + \pi_{1}\left( \mathbf{x} \right) + \pi_{2}\left( \mathbf{x} \right) = 1.
\tag{3.13}
\end{equation}

\textbf{Proof.} Each formula is (3.7) with (3.8)--(3.9), integrated over the boundary coordinate: on the arcs \(\int_{\varphi_{1}}^{\varphi_{2}}\sin\left( \nu_{m}\left( \psi - \varphi_{1} \right) \right)\, d\psi = \nu_{m}^{- 1}\left( 1 - ( - 1)^{m} \right)\), supported on odd \(m\); on the edges the radial integral produces \(\mathcal{J}_{m,p}\) of (3.2), with all \(m\). Substituting \(N_{m,p} = \frac{1}{2}\Delta\varphi\,\eta_{m,p}\) gives the displayed constants. Replacing \(e^{- \mu t}\) by \(\mu^{- 1}\) is the time integral \(\int_{0}^{\infty}\), valid since \(\tau < \infty\) a.s. The identity (3.13) is the statement that \(\sum_{\bullet}^{}\pi_{\bullet} = \sum_{m\ odd,p}^{}c_{m,p}\phi_{m,p} = 1\); it follows from the radial Lagrange identity \(\lbrack re^{- \alpha r^{2}}\mathcal{R}_{m,p}'\rbrack_{r_{1}}^{r_{2}} = \nu_{m}^{2}\,\mathcal{J}_{m,p} - \frac{2\mu_{m,p}}{\sigma^{2}}\mathcal{K}_{m,p}\), obtained by integrating the Sturm--Liouville form of (2.15), which converts the arc and edge contributions of (3.10)--(3.12) into \(c_{m,p}\) exactly. \(\quad\quad \ensuremath{\square}\)

\textbf{Remark 4 (which edge is favored: even versus odd harmonics).} Comparison of the two formulae in (3.12) shows that the edge densities differ only through the parity factor \(( - 1)^{m + 1}\). Consequently the \emph{total} edge flux carries only odd harmonics,

\[f_{\tau}^{\, 1} + f_{\tau}^{\, 2} = \frac{2\sigma^{2}}{\Delta\varphi}\sum_{m\ odd}^{}{\sum_{p}^{}\frac{\nu_{m}\mathcal{J}_{m,p}}{\eta_{m,p}}}e^{- \mu_{m,p}t}\phi_{m,p}\left( \mathbf{x} \right),\]

consistently with the absence of even modes from \(\mathbb{P}_{\mathbf{x}}(\tau < \infty) = 1\), whereas the \emph{difference} --- the imbalance between the two edges, i.e.~the answer to ``through which radial edge is the sector more likely to be left'' --- is carried entirely by the even harmonics,

\[f_{\tau}^{\, 1} - f_{\tau}^{\, 2} = \frac{2\sigma^{2}}{\Delta\varphi}\sum_{m\ even}^{}{\sum_{p}^{}\frac{\nu_{m}\mathcal{J}_{m,p}}{\eta_{m,p}}}e^{- \mu_{m,p}t}\phi_{m,p}\left( \mathbf{x} \right).\]

The even angular modes, invisible to the survival probability, encode precisely the directional asymmetry of the exit. They are activated by an off-center starting point \(\varphi \neq \frac{1}{2}\left( \varphi_{1} + \varphi_{2} \right)\), through the factor \(\sin\left( \nu_{m}\left( \varphi - \varphi_{1} \right) \right)\) in \(\phi_{m,p}\left( \mathbf{x} \right)\), and vanish on the angular bisector, where the two edges are hit with equal probability by symmetry.

Figure 2 illustrates the boundary-resolved exit-time densities of Corollary 4, the boundary-hitting probabilities and the survival probability, with a comparison to Monte Carlo simulation.

\begin{figure}[htbp]
\centering
\includegraphics[width=\textwidth]{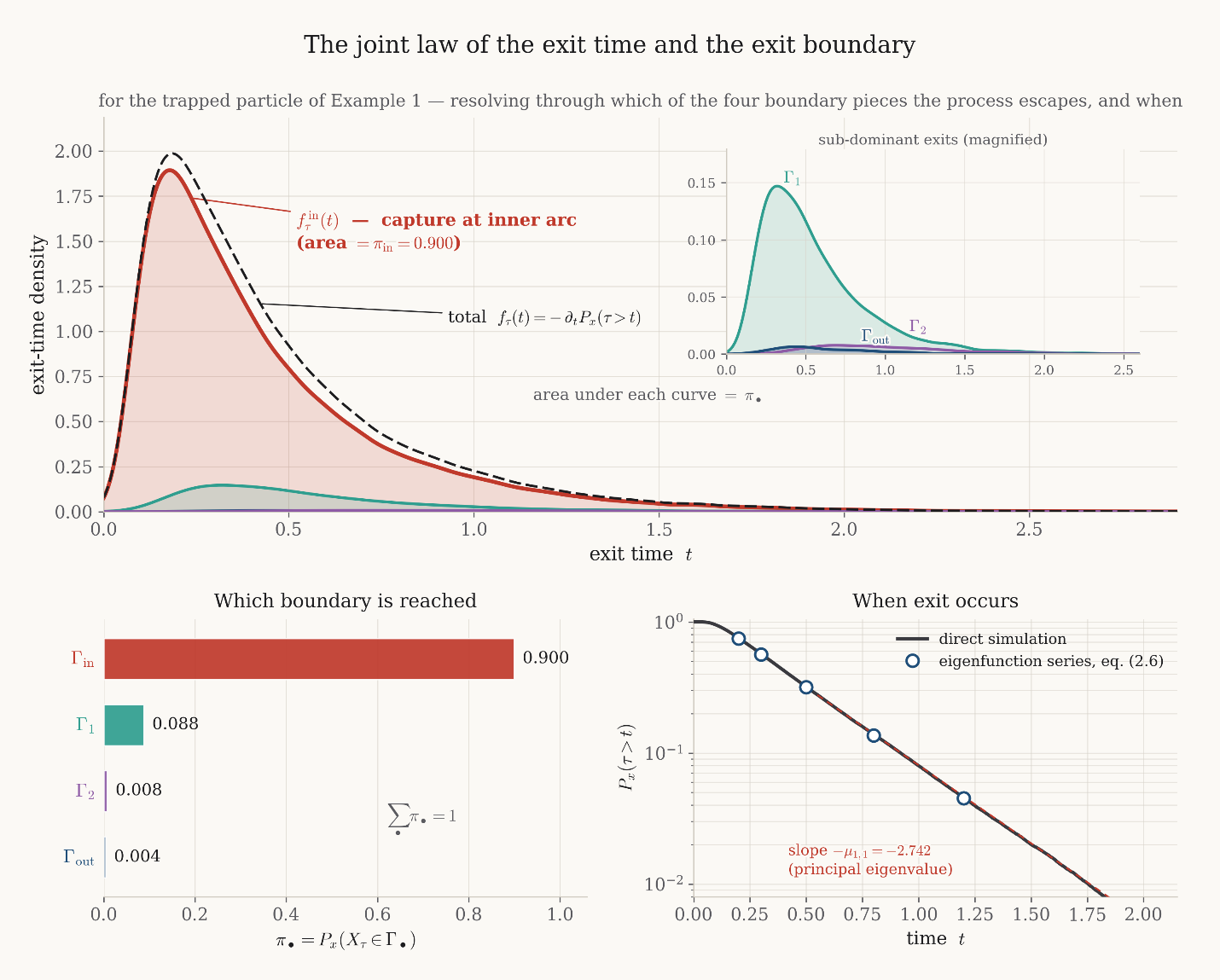}
\caption{\emph{The joint law of the exit time} \(\tau\) \emph{and the exit boundary} \(X_{\tau}\) \emph{for the trapped particle of the first example.} Boundaries are color-coded as in Figure 1. \textbf{(a)} The four boundary-resolved exit-time sub-densities \(f_{\tau}^{\bullet}(t)\) of Corollary 4: the area beneath each equals the corresponding hitting probability \(\pi_{\bullet}\), and their sum is the total exit-time density \(- \partial_{t}P_{x}(\tau > t)\) (dashed). The inset magnifies the three sub-dominant boundaries. \textbf{(b)} The boundary-hitting probabilities \(\pi_{\bullet} = P_{x}\left( X_{\tau} \in \Gamma_{\bullet} \right)\) from (3.10)--(3.13) --- here \(\pi_{in} = 0.900\), \(\pi_{1} = 0.088\), \(\pi_{2} = 0.008\), \(\pi_{out} = 0.004\) --- summing to unity. \textbf{(c)} The survival probability \(P_{x}(\tau > t)\) on a logarithmic scale: the solid curve is a direct Monte Carlo estimate from \(10^{5}\) trajectories, the open circles are the eigenfunction series (2.6) reported in Table 2, and the dashed guide carries the asymptotic slope \(- \mu_{1,1} = - 2.742\) fixed by the principal eigenvalue. The sub-densities in (a) and the probabilities in (b) are estimated from the same \(10^{5}\) simulated trajectories.}\label{fig2}
\end{figure}

\textbf{Remark 5 (quasi-stationarity).} By (3.4), the law of \(X_{t}\) conditioned on \(\{\tau > t\}\) converges as \(t \rightarrow \infty\) to the quasi-stationary distribution proportional to \(\phi_{1,1}\left( \mathbf{x} \right)\,\rho\left( \mathbf{x} \right)\) on \(D\), with mean residual lifetime \(\mu_{1,1}^{- 1}\); \(\mu_{1,1}\) is the principal Dirichlet eigenvalue of \(\mathcal{- L}\) on the annular sector, obtained from (2.9) at \(m = 1\).

\section{The reversible correlated Ornstein--Uhlenbeck process on a principal-axis rectangle}\label{sec:4}

We now treat a genuinely correlated planar OU process and establish the companion to Proposition 1. As anticipated in the introduction, curvature of the boundary is incompatible with correlation: the curved-domain construction of Section 2 separates only for the isotropic generator. The complementary tractable situation is a rectangle aligned with the principal axes of the process, on which the special functions are parabolic cylinder functions. The survival probability there factorizes into one-dimensional exit problems and so carries no new information; the genuinely two-dimensional object --- and the one we have not seen tabulated in this form --- is the joint law of the exit time and the exit \emph{side}, which does not factorize. The rectangle case is included not as a second application but to exhibit precisely what survives of the annular-sector method once isotropy is abandoned: separability persists only after passing to the principal axes, and the boundary-resolved law remains non-factorizing.

Let \(X = \left( X_{t} \right)_{t \geq 0}\) be the planar OU process

\begin{equation}
dX_{t} = - B\, X_{t}\, dt + \Sigma\, dW_{t},\quad\quad X_{t} \in \mathbb{R}^{2},\quad A: = \Sigma\Sigma^{\top} \succ 0,
\tag{4.1}
\end{equation}

with \(B\) a \(2 \times 2\) drift matrix having eigenvalues of positive real part (mean reversion), \(W\) a two-dimensional standard Brownian motion, and generator

\begin{equation}
\mathcal{L =}\frac{1}{2}\,\nabla \cdot (A\nabla) - (Bx) \cdot \nabla.
\tag{4.2}
\end{equation}

\subsection{Decorrelation and reduction to independent motions}\label{sec:4.1}

The eigenfunction method requires \(\mathcal{L}\) to be self-adjoint in a weighted \(L^{2}\) space, i.e.~\(X\) to be reversible. For a linear diffusion this is a condition relating drift and diffusion.

\textbf{Proposition 3 (decorrelation).} \emph{The process (4.1) is reversible with respect to its Gaussian invariant law if and only if}

\begin{equation}
BA = AB^{\top},
\tag{4.3}
\end{equation}

\emph{in which case} \(B\) \emph{is self-adjoint for the inner product} \(\langle u,w\rangle_{A^{- 1}}: = u^{\top}A^{- 1}w\)\emph{, equivalently} \(A^{- 1}B = B^{\top}A^{- 1}\)\emph{, and hence possesses real positive eigenvalues} \(\lambda_{1},\lambda_{2}\) \emph{with eigenvectors} \(u_{1},u_{2}\)\emph{, normalized to} \(\left| u_{i} \right| = 1\)\emph{, that are orthogonal in that inner product. Setting}

\begin{equation}
\mathbf{y} = M\mathbf{x},\quad\quad M = \lbrack\, u_{1}\ \ u_{2}\,\rbrack^{- 1},\quad\quad s_{i}^{2} = \frac{1}{u_{i}^{\top}A^{- 1}u_{i}},
\tag{4.4}
\end{equation}

\emph{transforms} \(X\) \emph{into a process} \(Y = MX\) \emph{whose two components are independent one-dimensional OU processes with distinct mean-reversion rates and volatilities,}

\begin{equation}
dY_{i,t} = - \lambda_{i}\, Y_{i,t}\, dt + s_{i}\, d{\widetilde{W}}_{i,t},\quad\quad i = 1,2,\quad{\widetilde{W}}_{1}\bot{\widetilde{W}}_{2},
\tag{4.5}
\end{equation}

\emph{and the generator into the decoupled sum}

\begin{equation}
\widetilde{\mathcal{L}} = \sum_{i = 1}^{2}\left( \frac{s_{i}^{2}}{2}\,\partial_{y_{i}y_{i}} - \lambda_{i}\, y_{i}\,\partial_{y_{i}} \right).
\tag{4.6}
\end{equation}

\emph{The squared volatilities} \(s_{i}^{2}\) \emph{in (4.4) are the diagonal entries of} \(MAM^{\top}\)\emph{; the rates} \(\lambda_{i}\) \emph{are intrinsic to} \(B\)\emph{; the numerical values of the} \(s_{i}\) \emph{and of the rectangle endpoints (4.7) depend on the coordinate normalization chosen for the principal variables, but the exit laws on the original parallelogram are invariant under this reparameterization.}

\textbf{Proof.} Reversibility of (4.1) is self-adjointness of \(\mathcal{L}\) in \(L^{2}(\pi)\), with \(\pi \propto e^{- \frac{1}{2}x^{\top}V^{- 1}x}\) the Gaussian invariant density and \(V\) the stationary covariance solving the Lyapunov equation \(BV + VB^{\top} = A\). The process is reversible if and only if \(BV\) is symmetric, equivalently \(BV = \frac{1}{2}\left( BV + VB^{\top} \right) = \frac{1}{2}A\), whence \(V = \frac{1}{2}B^{- 1}A\); symmetry of \(V\) is in turn equivalent to \(BA = AB^{\top}\), which is (4.3). Under (4.3) the matrix \(A^{- 1}B\) is symmetric, so \(B\) is self-adjoint for \(\langle \cdot , \cdot \rangle_{A^{- 1}}\) and admits real eigenvalues \(\lambda_{1},\lambda_{2}\) --- positive, since \(B\) generates a mean-reverting process --- with \(A^{- 1}\)-orthogonal eigenvectors \(u_{1},u_{2}\). Writing \(P = \left\lbrack \, u_{1}\ u_{2}\, \right\rbrack\) and \(M = P^{- 1}\), the similarity \(MBM^{- 1} = \operatorname{diag}\left( \lambda_{1},\lambda_{2} \right)\) diagonalizes the drift, while \(\left( MAM^{\top} \right)^{- 1} = P^{\top}A^{- 1}P\) has \((i,j)\) entry \(u_{i}^{\top}A^{- 1}u_{j}\), which is diagonal by orthogonality; hence \(MAM^{\top} = \operatorname{diag}\left( s_{1}^{2},s_{2}^{2} \right)\) with \(s_{i}^{2} = \left( u_{i}^{\top}A^{- 1}u_{i} \right)^{- 1}\). Thus \(Y = MX\) has both a diagonal drift and a diagonal diffusion, i.e.~its components are the independent one-dimensional OU processes (4.5) with decoupled generator (4.6). \(\quad\quad \ensuremath{\square}\)

The map \(M\) diagonalizes the drift and the diffusion simultaneously. Two consequences must be emphasized. First, although the transformed dynamics (4.5) are independent, the observable process \(X = M^{- 1}Y\) and, with it, the correlation structure \((B,A)\) are encoded in \(M\); the decorrelation is a change of coordinates, not a change of model. Second, the natural domain on which (4.6) admits a product solution is a rectangle in the \(\mathbf{y}\)-variables,

\begin{equation}
\mathcal{D = \{}\left( y_{1},y_{2} \right):\ \ell_{1} \leq y_{1} \leq h_{1},\ \ell_{2} \leq y_{2} \leq h_{2}\},
\tag{4.7}
\end{equation}

whose preimage \(\mathcal{P}: = M^{- 1}\mathcal{D}\) in the original coordinates is a parallelogram. Thus, within this self-adjoint, separable eigenfunction framework, just as the curved boundary of Section 2 required isotropy, a reversibly correlated process separates naturally on a rectangle aligned with its principal axes --- a parallelogram in the observed coordinates. We write \(\tau = \inf\{ t \geq 0:X_{t}\notin \mathcal{P}\} = \inf\{ t \geq 0:Y_{t}\notin \mathcal{D}\}\).

\subsection{The one-dimensional building block}\label{sec:4.2}

Each coordinate of (4.5) is a one-dimensional OU process on an interval, with its own rate \(\lambda\) and volatility \(s\), reversible with respect to \(\rho_{\lambda,s}(y) \propto e^{- \lambda y^{2}/s^{2}}\); its two-sided exit eigenproblem is classical \citep{alili2005,dinardo2001,linetsky2004}. We record it in the present notation, with the volatility retained explicitly.

\textbf{Lemma 1.} \emph{Let} \(\mathcal{L}_{\lambda,s} = \frac{s^{2}}{2}\partial_{yy} - \lambda y\partial_{y}\) \emph{on} \(\left( \ell,h \right)\) \emph{with Dirichlet conditions,} \(\lambda,s > 0\)\emph{. Its eigenpairs} \((\theta_{k},\chi_{k})_{k \geq 1}\) \emph{are}

\begin{equation}
\theta_{k} = \lambda\,\zeta_{k},\quad\quad\chi_{k}(y) = H_{\zeta_{k}}\left( \frac{\sqrt{\lambda}}{s}y \right)\, H_{\zeta_{k}}\left( - \frac{\sqrt{\lambda}}{s}\ell \right) - H_{\zeta_{k}}\left( - \frac{\sqrt{\lambda}}{s}y \right)\, H_{\zeta_{k}}\left( \frac{\sqrt{\lambda}}{s}\ell \right),
\tag{4.8}
\end{equation}

\emph{where} \(H_{\zeta}\) \emph{is the Hermite function of order} \(\zeta\) \emph{and the orders} \(\zeta_{k}\) \emph{are the increasing roots of the two-point equation}

\begin{equation}
H_{\zeta}\left( \frac{\sqrt{\lambda}}{s}\ell \right)\, H_{\zeta}\left( - \frac{\sqrt{\lambda}}{s}h \right) - H_{\zeta}\left( - \frac{\sqrt{\lambda}}{s}\ell \right)\, H_{\zeta}\left( \frac{\sqrt{\lambda}}{s}h \right) = 0.
\tag{4.9}
\end{equation}

\emph{Equivalently, through} \(H_{\zeta}(x) = 2^{\zeta/2}e^{x^{2}/2}D_{\zeta}\left( x\sqrt{2} \right)\)\emph{, the eigenfunctions are parabolic cylinder (Weber) functions} \(D_{\zeta}\) \emph{of non-integral order, with each argument} \(\frac{\sqrt{\lambda}}{s}( \cdot )\) \emph{replaced by} \(\frac{\sqrt{2\lambda}}{s}( \cdot )\)\emph{. The family} \(\{\chi_{k}\}\) \emph{is complete and orthogonal in} \(L^{2}\left( \left( \ell,h \right),\rho_{\lambda,s} \right)\)\emph{.}

\textbf{Proof.} With \(\xi = \sqrt{\lambda}\, y/s\) the eigen-ODE \(\frac{s^{2}}{2}\chi'' - \lambda y\chi' + \theta\chi = 0\) becomes Hermite's equation \(\chi_{\xi\xi} - 2\xi\chi_{\xi} + 2(\theta/\lambda)\chi = 0\), invariant under \(\xi \mapsto - \xi\); for non-integral order \(\theta/\lambda\) its two independent solutions are \(H_{\theta/\lambda}(\xi)\) and \(H_{\theta/\lambda}( - \xi)\). The combination (4.8) vanishes at \(y= \ell\) by construction and at \(y = h\) precisely when (4.9) holds. Regularity of both endpoints (\(\ell,h\) finite) gives, by Sturm--Liouville theory, a real simple spectrum increasing to \(+ \infty\) and a complete orthogonal eigenbasis in \(L^{2}\left( \rho_{\lambda,s} \right)\). \(\quad\quad \ensuremath{\square}\)

We attach a superscript \((i)\) to the quantities of Lemma 1 for the interval \(\left( \ell_{i},h_{i} \right)\), rate \(\lambda_{i}\) and volatility \(s_{i}\), so that \(\theta_{i,k} = \lambda_{i}\zeta_{i,k}\) and \(\rho_{i} \propto e^{- \lambda_{i}y^{2}/s_{i}^{2}}\), and we write the one-dimensional projection coefficients, survival functions and boundary fluxes

\begin{equation}
c_{k}^{(i)} = \frac{\langle 1,\chi_{k}^{(i)}\rangle_{\rho_{i}}}{\parallel \chi_{k}^{(i)} \parallel_{\rho_{i}}^{2}},\quad\quad S^{(i)}(t;y) = \sum_{k \geq 1}^{}c_{k}^{(i)}\, e^{- \theta_{i,k}t}\,\chi_{k}^{(i)}(y)
\tag{4.10}
\end{equation}

\begin{equation}
\begin{gathered}
f_{h_{i}}^{(i)}(t;y) = - \frac{s_{i}^{2}}{2}\,\rho_{i}\left( h_{i} \right)\sum_{k \geq 1}^{}\frac{\left( \chi_{k}^{(i)} \right)'\left( h_{i} \right)}{\parallel \chi_{k}^{(i)} \parallel_{\rho_{i}}^{2}}\, e^{- \theta_{i,k}t}\,\chi_{k}^{(i)}(y) \\
f_{\ell_{i}}^{(i)}(t;y) = + \frac{s_{i}^{2}}{2}\,\rho_{i}\left( \ell_{i} \right)\sum_{k \geq 1}^{}\frac{\left( \chi_{k}^{(i)} \right)'\left( \ell_{i} \right)}{\parallel \chi_{k}^{(i)} \parallel_{\rho_{i}}^{2}}\, e^{- \theta_{i,k}t}\,\chi_{k}^{(i)}(y)
\end{gathered}
\tag{4.11}
\end{equation}

\(f_{h_{i}}^{(i)}\) and \(f_{\ell_{i}}^{(i)}\) being the densities of exit of the \(i\)-th coordinate through its right and left endpoint, respectively (both non-negative as one-dimensional exit densities; for the principal mode \(\left( \chi_{k}^{(i)} \right)'\left( h_{i} \right) < 0 < \left( \chi_{k}^{(i)} \right)'\left( \ell_{i} \right)\), which fixes the sign of the leading term). The diffusion coefficient \(s_{i}^{2}\) now appears explicitly in the flux prefactor, in place of the isotropic \(\sigma^{2}\) of Section 3. The norms \(\parallel \chi_{k}^{(i)} \parallel_{\rho_{i}}^{2}\) and the coefficients \(c_{k}^{(i)}\) have closed forms obtained from the Sturm--Liouville boundary-term device of Remark 1, now with Hermite functions in place of confluent hypergeometric functions.

\subsection{Survival probability and the joint law of exit time and exit side}\label{sec:4.3}

Because the generator (4.6) decouples, the two coordinate processes are independent, and the survival probability is a product.

\textbf{Proposition 4 (survival probability).}

\emph{For} \(\mathbf{y} = \left( y_{1},y_{2} \right)\in \mathcal{D}\)\emph{,}

\begin{equation}
\begin{gathered}
\mathbb{P}_{\mathbf{y}}(\tau > t) = S^{(1)}\left( t;y_{1} \right)\, S^{(2)}\left( t;y_{2} \right) \\
= \sum_{j \geq 1}^{}{\sum_{k \geq 1}^{}c_{j}^{(1)}}c_{k}^{(2)}\, e^{- \left( \theta_{1,j} + \theta_{2,k} \right)t}\,\chi_{j}^{(1)}\left( y_{1} \right)\,\chi_{k}^{(2)}\left( y_{2} \right)
\end{gathered}
\tag{4.12}
\end{equation}

Indeed \(\{\tau > t\} = \{\tau^{(1)} > t\} \cap \{\tau^{(2)} > t\}\) with \(\tau^{(i)}\) the exit time of the independent coordinate \(i\). The double sum is a product of two one-dimensional spectral series; in this sense the survival probability of the correlated process is reducible to known one-dimensional objects, and the same is true of the exit-time density \(- \partial_{t}\) and the moments. The genuinely two-dimensional content lies elsewhere, in the resolution of \emph{which side} of the rectangle is reached, since the exit through a given edge requires the \emph{other} coordinate to have survived up to the same random instant --- a coupling that does not factorize.

Denote the four open edges of \(\mathcal{D}\) by \(E_{1}^{+} = \{ y_{1} = h_{1}\}\), \(E_{1}^{-} = \{ y_{1} = \ell_{1}\}\), \(E_{2}^{+} = \{ y_{2} = h_{2}\}\), \(E_{2}^{-} = \{ y_{2} = \ell_{2}\}\) (each with the complementary coordinate ranging over its open interval). Exit at a corner has probability zero, so these four events partition \(\{ X_{\tau} \in \partial\mathcal{P}\}\) almost surely.

\textbf{Proposition 5 (joint law of exit time and exit side).} \emph{For} \(\mathbf{y}\in \mathcal{D}\)\emph{, the four sub-densities of the exit time} \(\tau\) \emph{on the edges of} \(\mathcal{D}\) \emph{are products of a one-dimensional exit flux of the coordinate that escapes and the one-dimensional survival function of the coordinate that remains,}

\begin{equation}
f_{E_{1}^{\pm}}\left( t;\mathbf{y} \right) = f_{h_{1}/\ell_{1}}^{(1)}\left( t;y_{1} \right)\, S^{(2)}\left( t;y_{2} \right),\quad\quad f_{E_{2}^{\pm}}\left( t;\mathbf{y} \right) = S^{(1)}\left( t;y_{1} \right)\, f_{h_{2}/\ell_{2}}^{(2)}\left( t;y_{2} \right),
\tag{4.13}
\end{equation}

\emph{and their sum equals} \(- \partial_{t}\mathbb{P}_{\mathbf{y}}(\tau > t)\)\emph{. More finely, the joint density of the pair (exit time} \(t\)\emph{, position} \(y_{2}'\) \emph{along the edge) on} \(E_{1}^{+}\) \emph{is} \(f_{h_{1}}^{(1)}\left( t;y_{1} \right)\, q^{(2)}\left( t;y_{2},y_{2}' \right)\)\emph{, with} \(q^{(2)}\) \emph{the killed transition density of coordinate} \(2\) \emph{with respect to Lebesgue measure, and analogously on the other edges.}

\emph{Consequently the probabilities of leaving through each edge are}

\begin{equation}
\pi_{E_{1}^{+}}\left( \mathbf{y} \right) = - \frac{s_{1}^{2}}{2}\,\rho_{1}\left( h_{1} \right)\sum_{j \geq 1}^{}{\sum_{k \geq 1}^{}\frac{\left( \chi_{j}^{(1)} \right)'\left( h_{1} \right)\, c_{k}^{(2)}}{\parallel \chi_{j}^{(1)} \parallel_{\rho_{1}}^{2}}}\,\frac{\chi_{j}^{(1)}\left( y_{1} \right)\,\chi_{k}^{(2)}\left( y_{2} \right)}{\lambda_{1}\zeta_{1,j} + \lambda_{2}\zeta_{2,k}}
\tag{4.14}
\end{equation}

\emph{and, with the obvious changes of endpoint and coordinate,} \(\pi_{E_{1}^{-}},\pi_{E_{2}^{+}},\pi_{E_{2}^{-}}\)\emph{; they satisfy}

\begin{equation}
\pi_{E_{1}^{+}}\left( \mathbf{y} \right) + \pi_{E_{1}^{-}}\left( \mathbf{y} \right) + \pi_{E_{2}^{+}}\left( \mathbf{y} \right) + \pi_{E_{2}^{-}}\left( \mathbf{y} \right) = 1.
\tag{4.15}
\end{equation}

\textbf{Proof.} By independence of the two coordinates, the event of leaving \(\mathcal{D}\) through \(E_{1}^{+}\) in \(\lbrack t,t + dt)\) with the second coordinate near \(y_{2}'\) is the intersection of \(\{\)coordinate \(1\) exits at \(h_{1}\) in \(\lbrack t,t + dt)\}\), of density \(f_{h_{1}}^{(1)}\left( t;y_{1} \right)\), and \(\{\)coordinate \(2\) is near \(y_{2}'\) and has not exited\(\}\), of density \(q^{(2)}\left( t;y_{2},y_{2}' \right)\); the product is the stated joint density. Integrating \(y_{2}'\) over \(\left( \ell_{2},h_{2} \right)\) replaces \(q^{(2)}\) by \(S^{(2)}\) and gives (4.13). Summing the four edge densities and using \(- \partial_{t}S^{(i)} = f_{h_{i}}^{(i)} + f_{\ell_{i}}^{(i)}\) together with the product rule applied to (4.12) yields \(- \partial_{t}\mathbb{P}_{\mathbf{y}}(\tau > t)\). Integrating (4.13) over \(t \in (0,\infty)\), substituting the spectral forms (4.10)--(4.11) and using \(\int_{0}^{\infty}e^{- \left( \theta_{1,j} + \theta_{2,k} \right)t}\, dt = \left( \theta_{1,j} + \theta_{2,k} \right)^{- 1} = \left( \lambda_{1}\zeta_{1,j} + \lambda_{2}\zeta_{2,k} \right)^{- 1}\) gives (4.14). Identity (4.15) is the statement that exit is almost sure, \(\int_{0}^{\infty}\left( - \partial_{t}\mathbb{P}_{\mathbf{y}}(\tau > t) \right)\, dt = 1\). \(\quad\quad \ensuremath{\square}\)

The factor \(\left( \lambda_{1}\zeta_{1,j} + \lambda_{2}\zeta_{2,k} \right)^{- 1}\) in (4.14) is the signature of the coupling: the side probabilities are \emph{not} products of one-dimensional quantities, the two spectra entering inseparably through the sum of eigenvalues, even though the survival probability (4.12) is a clean product. This is the correlated-process counterpart of the boundary decomposition of Section 3, and it is here that the analytically non-trivial, genuinely two-dimensional content of the correlated problem resides.

\textbf{Remark 6 (the geometry/operator dichotomy).} Propositions 1 and 5 delimit the two canonical separable geometries for planar OU exit. The curved annular sector demands isotropy and yields confluent hypergeometric (Whittaker) radial functions with a non-trivially coupled, non-factorizing survival probability; the principal-axis rectangle admits arbitrary reversible correlation and yields parabolic cylinder functions with a factorizing survival probability but a non-factorizing side law. The irreversible case \(BA \neq AB^{\top}\) --- in particular a rotational (spiralling) drift, for which \(B\) has complex eigenvalues --- lies outside both frameworks and is the genuine frontier flagged in the introduction.

Figure 3 summarizes this geometry/operator dichotomy, contrasting the isotropic curved problem with the correlated rectilinear one.

\begin{figure}[htbp]
\centering
\includegraphics[width=\textwidth]{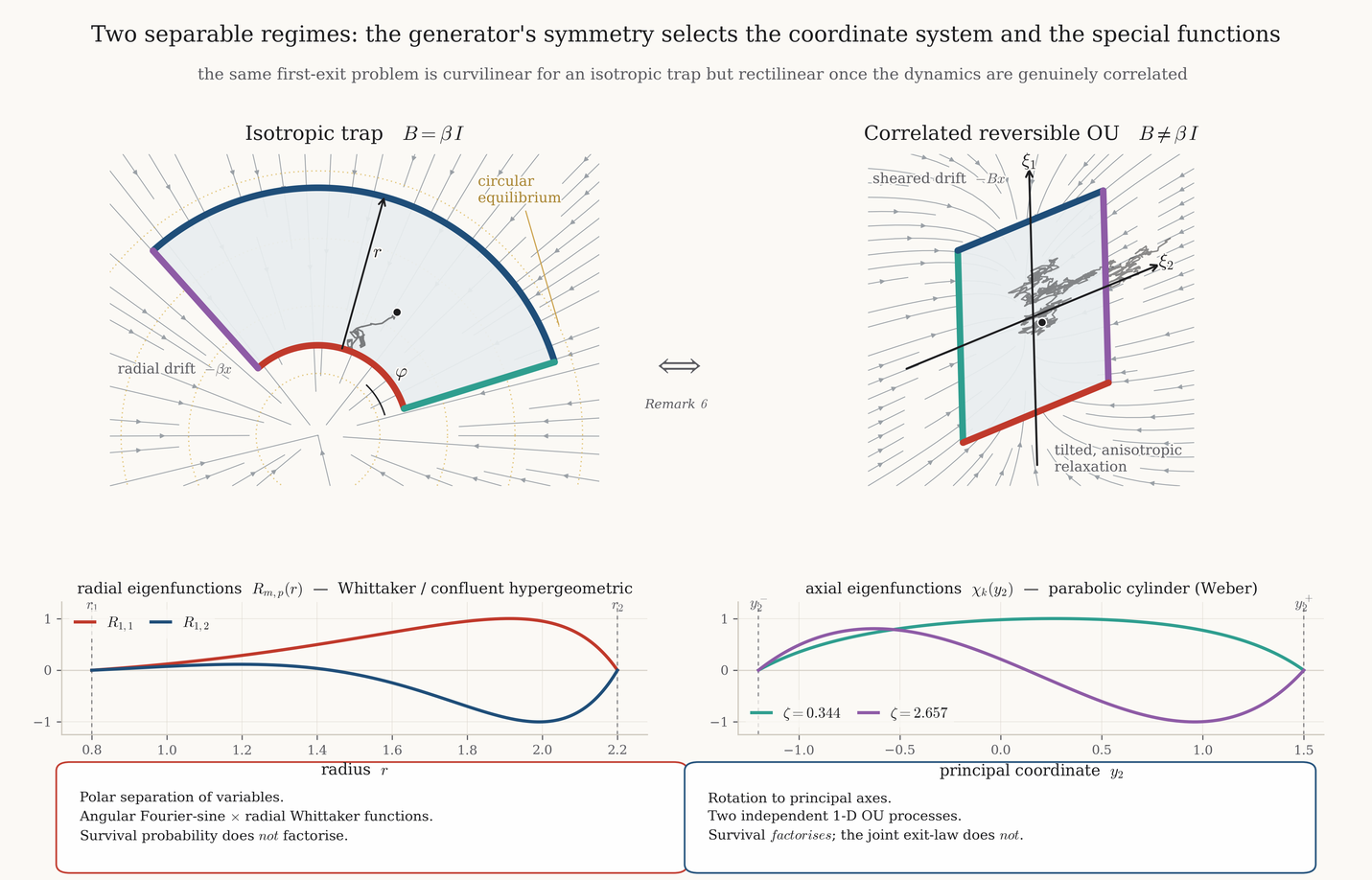}
\caption{\emph{The geometry/operator dichotomy for planar Ornstein--Uhlenbeck exit (Remark 6).} \textbf{Left:} the isotropic trap \(B = \beta I\). The drift \(- \beta x\) is radial and the equilibrium is circular, so the generator separates in polar coordinates \((r,\varphi)\) and the natural domain is the annular sector. \textbf{Right:} a genuinely correlated reversible process \(B \neq \beta I\). The drift \(- Bx\) is sheared and the generator separates along the straight principal axes \(\xi_{1},\xi_{2}\) (Proposition 3), whose natural domain is the principal-axis rectangle --- a parallelogram in the observed coordinates. Streamlines indicate the drift, with a sample path overlaid in each panel. \textbf{Bottom:} the radial eigenfunctions \(R_{m,p}(r)\) of the curved problem are Whittaker (confluent hypergeometric) functions vanishing at \(r_{1}\) and \(r_{2}\) (shown: \(R_{1,1}\) and \(R_{1,2}\), with \(\mu_{1,1} \approx 2.74\), \(\mu_{1,2} \approx 6.70\)), whereas the axial eigenfunctions \(\chi_{k}\left( y_{2} \right)\) of the rectilinear problem are parabolic cylinder (Weber) functions of non-integral order (shown: \(\zeta = 0.344,\ 2.657\)). In the isotropic case the survival probability does not factorize (Proposition 1); in the correlated case it factorizes into one-dimensional factors (Proposition 4), yet the joint exit-side law does not (Proposition 5). Left panel: the sector of the first example; right panel and the orders \(\zeta\): the reversible correlated example of \S5.4.}\label{fig3}
\end{figure}

\section{Numerical implementation and validation}\label{sec:5}

The series of Propositions 1, 2 and 5 are explicit, but they involve special functions of non-integral order and eigenvalues defined implicitly; this section describes their numerical evaluation and validates the results against direct Monte Carlo simulation. We treat the annular-sector process of Sections 2--3 in detail; the rectangle of Section 4, assembled from one-dimensional parabolic-cylinder problems, is handled analogously and more simply.

\subsection{Computation of the spectrum and the series}\label{sec:5.1}

\emph{Special functions.} The Kummer and Tricomi functions \(M(a,b,z)\) and \(U(a,b,z)\) are evaluated in extended precision --- thirty to fifty significant digits, adaptively refined --- from their standard convergent and asymptotic series \citep[\S13]{mathews2022,olver2010}; the digits tabulated below are produced at this precision and then rounded. Extended precision is not a convenience but a necessity: for the higher radial modes \(a\) is large and negative, and the cancellation between the two solutions entering (2.9) --- the same near-boundary sensitivity discussed in Section 5.2 --- costs several digits in ordinary double precision. For an independent guarantee each value may be enclosed rigorously by interval arithmetic, or read directly off the shooting integration described below.

\emph{Eigenvalues and root-finding.} For each angular index \(m\), \(\mu_{m,p}\) is the p-th eigenvalue of the radial Dirichlet Sturm--Liouville problem; equivalently the corresponding parameter \(a\) is a non-spurious root of the two-radius determinant (2.9), with \(\mu_{m,p} = \beta\left( \nu_{m} - 2a_{m,p} \right)\). It is essential to root-find not on (2.9) itself but on its deflated form. Writing the second Frobenius solution at the origin as \(\widetilde{M}(a,b,z) = z^{1 - b}\, M(a - b + 1,2 - b,z)\) and using the connection formula \(U = \frac{\Gamma(1 - b)}{\Gamma(a - b + 1)}\, M + \frac{\Gamma(b - 1)}{\Gamma(a)}\,\widetilde{M}\), the determinant factorizes as

\[(2.9) = \frac{\Gamma(b - 1)}{\Gamma(a)} \cdot \widetilde{\Delta}(a),\,\,\,\widetilde{\Delta}(a): = M\left( a,b,z_{1} \right)\,\widetilde{M}\left( a,b,z_{2} \right) - M\left( a,b,z_{2} \right)\,\widetilde{M}\left( a,b,z_{1} \right),\]

in which the entire functions \(M\) and \(\widetilde{M}\) are linearly independent for the non-integral order \(b = \nu_{m} + 1\), so that \(\widetilde{\Delta}\) vanishes precisely at the genuine eigenvalues. When \(b = \nu_{m} + 1\) is an integer, the displayed Frobenius and connection representation is understood by analytic continuation in the angular order, equivalently by replacing the second local solution with the logarithmic limiting solution; the deflated determinant is then the corresponding analytic continuation, and the Sturm--Liouville spectrum is unchanged. The spurious zeros of (2.9) at \(a \in \left\{ 0, - 1, - 2,\ldots \right\}\) are exactly the zeros of the prefactor \(1/\Gamma(a)\) and are absent from \(\widetilde{\Delta}\); root-finding on \(\widetilde{\Delta}\) therefore requires no exclusion of integer neighborhoods --- a real advantage, since a genuine root may lie arbitrarily close to such a point (in the example of Section 5.2, \(a_{1,4} \approx - 9.97\), within 0.03 of the spurious value \(- 10\)). Since \(\mu_{m,p} > 0\) forces \(a < \nu_{m}/2\), we scan \(a\) downward from \(\nu_{m}/2 - \delta\) to \(a_{\min} = \left( \nu_{m} - \Lambda/\beta \right)/2\), with \(\Lambda\) the radial cutoff below, bracket the sign changes of \(\widetilde{\Delta}\) on a grid coarse enough to separate consecutive eigenvalues (a step of 0.1 suffices here, the gaps in \(a\) of about two or more in the range scanned, and widening with \(p\)), and refine each bracket by Brent's method to \(10^{- 12}\). Every root is cross-checked against the shooting integration of (2.15) from \(R\left( r_{1} \right) = 0\), \(R^{'}\left( r_{1} \right) = 1\): the index \(p\) is the number of interior zeros of the resulting solution, a monotone step function of \(\mu\), so that a bisection in \(\mu\) recovers each eigenvalue without reference to the determinant --- and never encounters the integer degeneracies --- confirming the closed form to all displayed digits (Table 1).

\emph{Quadrature and radial functionals.} The functionals \(K_{m,p}\), \(\eta_{m,p}\), \(J_{m,p}\) of (3.1)--(3.2) and the boundary derivatives \(R_{m,p}^{'}\left( r_{1} \right)\), \(R_{m,p}^{'}\left( r_{2} \right)\) are computed from a single, unnormalized representation of \(R_{m,p}\), the coefficients following from \(c_{m,p} = 4K_{m,p}/\left( \nu_{m}\,\Delta\varphi\,\eta_{m,p} \right)\). The eigenfunctions are deliberately not normalized before differentiation: every reported quantity --- the joint density (3.7) and the boundary fluxes (3.8)--(3.12) --- is a ratio in which the arbitrary scale of \(\phi_{m,p}\) cancels, the numerator carrying \(\phi_{m,p}^{2}\) and the norm \(N_{m,p}\) carrying \(\left\| \phi_{m,p} \right\|^{2}\), so that an explicit normalization would only insert a further, ill-conditioned division. For conditioning, \(R_{m,p}\) is instead rescaled to order unity before any derivative or integral is formed, which leaves the ratios unchanged; the boundary slopes are taken analytically from \(\frac{dM}{dz} = (a/b)M(a + 1,b + 1,z)\), \(\frac{dU}{dz} = - a\, U(a + 1,b + 1,z)\) and \(\frac{dz}{dr} = 2\alpha r\), never by differencing a function that vanishes at the endpoint. The integrands are smooth on \(\left\lbrack r_{1},r_{2} \right\rbrack\) with at most \(p\) interior oscillations and are evaluated by Gauss--Legendre quadrature with the Gaussian factor folded into the integrand (\(n \approx 2p + 10\) nodes, the count doubled until the relative change falls below \(10^{- 10}\)), a tolerance well beneath both the reporting precision and the series-truncation error. When \(R_{m,p}\) is obtained by shooting, the functionals are appended to the radial system as quadratures \(\frac{dI}{dr} = R_{m,p}(r)\, e^{- \alpha r^{2}}\, r\) and produced by the same integrator; \(\eta_{m,p}\) is checked independently against the Sturm--Liouville boundary-term formula of Remark 1, and the moment integrals against the confluent-hypergeometric primitives of \citep[\S13.10]{olver2010}.

\emph{Truncation.} The double series are truncated at angular order \(m \leq M\) and radial eigenvalue \(\mu_{m,p} \leq \Lambda\). For the survival probability, density and moments the factor \(e^{- \mu_{m,p}\, t}\) governs the tail, and the cutoff \(\Lambda = t_{\min}^{- 1}\ln(1/\varepsilon)\) --- with \(\varepsilon\) the target relative accuracy and \(t_{\min}\) the smallest time evaluated --- secures it: at \(\varepsilon = 10^{- 8}\) and \(t_{\min} = 0.2\) this gives \(\Lambda \approx 92\), so the value \(\Lambda = 200\) used below is conservative. Within each \(m\) the \(p\)-sum is stopped at the first index for which \(\left| c_{m,p}\, e^{- \mu_{m,p}\, t}\,\phi_{m,p}(x) \right|\) falls below \(\varepsilon\) times the partial sum over two or three consecutive terms, the consecutive-term test guarding against an accidental node of \(\phi_{m,p}\). The boundary-hitting probabilities are the binding case: their time integral removes the exponential damping, after which the radial sum converges only algebraically and the angular sum conditionally, like \(1/m\). We therefore fix \(M\) and \(\Lambda\) from these quantities and monitor convergence through the exact mass identity (3.13), driving \(\left| \pi_{out} + \pi_{in} + \pi_{1} + \pi_{2} - 1 \right|\) below the target as \(M\) and \(\Lambda\) are raised; this is displayed for the boundary-hitting probabilities in Table 3 below, and the acceleration used to reach the converged values is described next (Algorithm 1).

\emph{Sequence acceleration.} The two arc fluxes \(\pi_{in}\) and \(\pi_{out}\) are obtained by collapsing the radial series first: for each angular mode the oscillatory radial partial sums of the arc contribution are extrapolated, while the rapidly convergent edge series giving \(\pi_{1}\) and \(\pi_{2}\) are summed as they stand, leaving a single coefficient \(A_{m}\) per mode. The residual angular series decays only like \(\frac{1}{m}\), the algebraic rate set by the re-entrant corners where the arcs meet the straight edges, so a second pass is applied to the partial sums \(S_{M} = \sum_{m = 1}^{M}A_{m}\). The even- and odd-indexed angular subsequences are extrapolated separately, which preserves the reflection symmetry between the two straight edges and keeps the parities clean in the sense of Remark 4. The recursion below is seeded with \(\varepsilon_{- 1}^{(M)} = 0\) and \(\varepsilon_{0}^{(M)} = S_{M}\); the odd-numbered columns are auxiliary and the accelerated estimate is read from the even columns. We report the deepest stable even column \(\varepsilon_{2q}^{(M)}\) reached at \(M = 48\) (twenty-four terms per parity, \(q = 11\)). The full procedure is collected in Algorithm 1, and the stability of the accelerated estimate as the truncation grows is documented in Table 4.

\[\varepsilon_{- 1}^{(M)} = 0,\ \ \ \varepsilon_{0}^{(M)} = S_{M},\ \ \ \varepsilon_{k + 1}^{(M)} = \varepsilon_{k - 1}^{(M + 1)} + \frac{1}{\varepsilon_{k}^{(M + 1)} - \varepsilon_{k}^{(M)}}\]

\medskip
\noindent\textbf{Algorithm 1. Wynn $\varepsilon$-acceleration of the boundary fluxes.}
{\small
\begin{verbatim}
Input:  angular truncation M, radial truncation Lambda.
Output: accelerated fluxes pi_in, pi_out, pi_1, pi_2.

for m = 1 to M:                      # one coefficient per angular mode
    arc modes:  A_in[m], A_out[m] <- WynnEven( radial partial sums )
    edge modes: A_1[m],  A_2[m]   <- direct sum over mu <= Lambda

for f in { in, out, 1, 2 }:          # accelerate the angular series
    split A_f[1..M] into odd-m and even-m subsequences
    pi_f <- WynnEven( odd-m sums ) + WynnEven( even-m sums )

return pi_in, pi_out, pi_1, pi_2

function WynnEven(S):                 # S[M] = partial sums, indexed by M
    eps(-1, M) = 0,  eps(0, M) = S[M]
    for k = 0, 1, 2, ... :
        eps(k+1, M) = eps(k-1, M+1) + 1 / ( eps(k, M+1) - eps(k, M) )
    return the deepest stable even column eps(2q, M)
\end{verbatim}
}
\medskip

\subsection{Validation against Monte Carlo simulation}\label{sec:5.2}

We take \(\beta = 1\), \(\sigma = 0.7\) (so \(\alpha = \beta/\sigma^{2} \approx 2.041\)), inner and outer radii \(r_{1} = 0.8\), \(r_{2} = 2.2\), and angular limits \(\varphi_{1} = 0.3\), \(\varphi_{2} = 2.3\), so that \(\Delta\varphi = 2\) and the orders \(\nu_{m} = m\pi/2\) are non-integral. The process starts at \(r_{0} = 1.30\), \(\varphi_{0} = 1.00\), off the angular bisector, so that the even angular modes are active. The trap center --- the origin --- lies \emph{inside} the inner radius, so the mean-reverting drift points on average through the inner arc; this will dominate the exit.

The radial eigenvalues from the closed-form determinant (2.9), with the spurious integer roots removed, and from the shooting integration of (2.15) coincide to the precision shown.

\begin{table}[htbp]
\caption{Lowest radial eigenvalues $\mu_{m,p}$; the closed-form condition (2.9) and the shooting method agree to all four digits.}\label{tab1}
\centering
\begin{tabular}{@{}ccccc@{}}
\toprule
$m$ & $\mu_{m,1}$ & $\mu_{m,2}$ & $\mu_{m,3}$ & $\mu_{m,4}$ \\
\midrule
1 & 2.7423 & 6.6960 & 12.8779 & 21.5135 \\
2 & 3.7066 & 7.6541 & 13.8689 & 22.5206 \\
3 & 5.1986 & 9.2616 & 15.5388 & 24.2137 \\
\botrule
\end{tabular}
\end{table}

The survival probability (2.6) is compared with the simulation in Table 2. The eigenfunction series and the Brownian-bridge--corrected Monte Carlo estimate agree throughout, each analytic value falling within the 95\% confidence interval of its simulated counterpart, over more than a decade of decay in the survival probability. The factor \(e^{- \mu_{m,p}t}\) damps the high modes, so for the times tabulated a moderate truncation already resolves the series to the displayed precision.

\begin{table}[htbp]
\caption{Survival probability $\mathbb{P}_{\mathbf{x}}(\tau > T)$ from the eigenfunction series and from $5 \times 10^{5}$ Brownian-bridge--corrected Monte Carlo trajectories ($\Delta t = 0.002$); 95\% confidence intervals shown.}\label{tab2}
\centering
\begin{tabular}{@{}lll@{}}
\toprule
$T$ & analytic, eq.\ (2.6) & Monte Carlo (95\% CI) \\
\midrule
0.2 & 0.7495 & $0.7477 \pm 0.0012$ \\
0.3 & 0.5655 & $0.5639 \pm 0.0014$ \\
0.5 & 0.3178 & $0.3176 \pm 0.0013$ \\
0.8 & 0.1366 & $0.1359 \pm 0.0009$ \\
1.2 & 0.0453 & $0.0448 \pm 0.0006$ \\
\botrule
\end{tabular}
\end{table}

The four boundary-hitting probabilities (3.10)--(3.13) converge far more slowly: as properties of the boundary they probe the expansion of the constant initial datum near \(\partial D\), where it converges only conditionally. Table 3 exhibits this directly, listing the four probabilities and their sum at three successive truncations \((M,\Lambda) = (24,120),(36,200),(48,300)\), together with the mass defect \(\left| \pi_{out} + \pi_{in} + \pi_{1} + \pi_{2} - 1 \right|\) which, by the identity (3.13), measures the truncation error directly, without reference to any exact value.

\begin{table}[htbp]
\caption{Boundary-hitting probabilities and mass defect for the first example. The first three rows are the raw double partial sums at the indicated truncations; the $\varepsilon$-accelerated row is the same series after acceleration of the angular sequence; the final row is the bridge-corrected Monte Carlo estimate ($5 \times 10^{5}$ paths, 95\% half-widths at most $8 \times 10^{-4}$).}\label{tab3}
\centering
\begin{tabular}{@{}lcccccc@{}}
\toprule
$(M,\Lambda)$ & $\pi_{in}$ & $\pi_{out}$ & $\pi_{1}$ & $\pi_{2}$ & $\Sigma$ & $|\Sigma - 1|$ \\
\midrule
$(24, 120)$ & 0.9536 & 0.0026 & 0.0585 & 0.0256 & 1.0402 & $4.0\times10^{-2}$ \\
$(36, 200)$ & 0.8946 & 0.0029 & 0.0954 & 0.0228 & 1.0155 & $1.6\times10^{-2}$ \\
$(48, 300)$ & 0.8864 & 0.0031 & 0.1212 & 0.0086 & 1.0193 & $1.9\times10^{-2}$ \\
$\varepsilon$-accelerated & 0.8999 & 0.0044 & 0.0879 & 0.0077 & 1.0000 & $5.7\times10^{-5}$ \\
Monte Carlo & 0.9000 & 0.0043 & 0.0877 & 0.0080 & 1.0000 & --- \\
\botrule
\end{tabular}
\end{table}

The raw sums converge slowly and non-monotonically, exactly as the unit mass they reconstruct would lead one to expect. That constant violates the homogeneous Dirichlet condition on all four sides of the sector, so both the angular expansion in \(m\) and the radial expansion in \(p\) behave like the Fourier series of a discontinuity: their coefficients decay only like \(1/m\) and \(1/p\), and the partial sums oscillate about the limit in the Gibbs manner rather than approaching it monotonically. The mass defect is accordingly a few percent, equivalently a few parts in \(10^{2}\) at these truncations and is not monotone in \((M,\Lambda)\) --- a faithful diagnostic of the slow summation, not of any error in the closed forms, whose eigenvalues and mode integrals are independently exact (Table 1 and Remark 1).

The fourth row applies the Wynn \emph{\ensuremath{\varepsilon}}-algorithm to the angular sequence of radially-converged partial sums --- the oscillatory radial arc-sums are themselves extrapolated first, the rapidly-convergent edge-sums taken as they stand, and the even and odd angular subsequences accelerated separately so that the reflection symmetry of the two straight edges is preserved. The extrapolation collapses the mass defect to \(6 \times 10^{- 5}\) and reproduces the independent bridge-corrected Monte Carlo estimates of all four probabilities (final row) to within \(5 \times 10^{- 4}\), comfortably inside the simulation's confidence interval. The accelerated values \(0.8999\), \(0.0044\), \(0.0879\), \(0.0077\) are therefore the converged probabilities; a directly-summed value at any single truncation has simply not yet converged for these conditionally-summing boundary fluxes, and it is the accelerated figures, in agreement with the simulation, that are taken as the converged values.

The accelerated estimates are stable under refinement of the truncation. Table 4 lists the \emph{\ensuremath{\varepsilon}}-accelerated fluxes computed from partial sums carried to successively larger orders. The dominant flux \(\pi_{in}\) is already fixed in its third decimal at the coarsest truncation and the remaining fluxes settle by the second, while the mass defect falls to \(5.7 \times 10^{- 5}\), more than two orders of magnitude below the defect of the corresponding raw partial sums in Table 3.

\begin{table}[htbp]
\caption{Stability of the accelerated boundary-hitting probabilities under refinement. Each row applies the $\varepsilon$-acceleration of Algorithm 1 to the partial sums carried to the stated orders $M$ and $\Lambda$; the four fluxes and their sum are reported as in Table 3. The mass defect falls by more than two orders of magnitude across the range, and the last row reproduces the accelerated row of Table 3 to all displayed digits.}\label{tab4}
\centering
\begin{tabular}{@{}lcccccc@{}}
\toprule
$(M,\Lambda)$ & $\pi_{in}$ & $\pi_{out}$ & $\pi_{1}$ & $\pi_{2}$ & $\Sigma$ & $|\Sigma - 1|$ \\
\midrule
$(24, 120)$ & 0.9020 & 0.0040 & 0.0861 & 0.0086 & 1.0007 & $7.0\times10^{-4}$ \\
$(36, 200)$ & 0.9005 & 0.0043 & 0.0876 & 0.0078 & 1.0002 & $2.0\times10^{-4}$ \\
$(48, 300)$ & 0.8999 & 0.0044 & 0.0879 & 0.0077 & 1.0000 & $5.7\times10^{-5}$ \\
\botrule
\end{tabular}
\end{table}

The physical picture of Section 3 is borne out. With the trap center inside the inner radius the mean-reverting drift points, on average, through the inner arc, and exit occurs overwhelmingly there --- close to ninety per cent --- and almost never through the outer arc, below one half of one per cent. The pronounced asymmetry between the two radial edges, \(\pi_{1}\, \gg \,\pi_{2}\) with the start lying nearer the edge \(\varphi_{1}\), is reproduced through the even angular modes as in Remark 4.

The Monte Carlo reference uses the exact Gaussian transition of the Ornstein--Uhlenbeck process over a step \(\Delta t\), with no Euler discretization of the dynamics,

\begin{equation}
X_{t + \Delta t} = e^{- \beta\Delta t}X_{t} + \sqrt{\frac{\sigma^{2}}{2\beta}\left( 1 - e^{- 2\beta\Delta t} \right)}Z,\, Z \sim N\left( 0,I_{2} \right)
\tag{5.1}
\end{equation}

with \(N = 5 \times 10^{5}\) trajectories and \(\Delta t = 0.002\). The first time a step lands outside \(D\) the path is recorded as having exited, through the boundary piece it has crossed. Discrete sampling alone, however, cannot see an excursion that leaves \(D\) and returns within a single step, and discarding those excursions biases every exit probability. The error is of order \(\sqrt{\Delta t}\), not \(\Delta t\): over one step the process is diffusive, so the chance of having crossed a nearby boundary scales with the step's standard deviation \(\sigma\sqrt{\Delta t}\) rather than its variance. We remove this leading bias by a Brownian-bridge correction. Conditional on the endpoints \(X_{t}\) and \(X_{t + \Delta t}\) of a step that stays inside \(D\), the probability that the connecting diffusion bridge has touched a locally straight boundary lying at perpendicular distances \(d\) and \(d^{'}\) from the two endpoints is \(\exp\left( - 2dd^{'}/\left( \sigma^{2}\Delta t \right) \right)\); we evaluate this for each arc --- the perpendicular distance being the radial gap --- and each edge --- the tangential gap \(r\sin\left( \varphi - \varphi_{edge} \right)\) --- and record the path as having exited through that piece. When more than one boundary piece carries appreciable crossing probability over the same step, the competing pieces are resolved by a single allocation: one uniform variate on the unit interval is drawn, and the interval is partitioned into consecutive subintervals whose lengths equal the individual bridge-crossing probabilities, the remaining length representing no crossing. The path is recorded as having exited through the piece whose subinterval contains the variate, and is otherwise continued into the interior of the sector. Because the pieces are separated by distances large compared with the typical one-step displacement, the crossing events are mutually exclusive up to a term exponentially small in the squared separation measured in units of the step variance; the subinterval lengths sum to less than one, the probability of a joint double crossing is of the same negligible order, and the measure-zero corner events are disregarded. The rule therefore coincides with independent Bernoulli absorption on each piece, with any double crossing assigned to the nearer piece. The correction is exact for a straight boundary and a Brownian bridge and accurate to higher order for the gently curved arcs and the small drift over a step; with it the discrete bias drops from \(O\left( \sqrt{\Delta t} \right)\) to \(O(\Delta t)\) and the simulated survival becomes flat in \(\Delta t\) at the value of the series. Every estimate is reported with its 95\% confidence interval \(\widehat{p} \pm 1.96\sqrt{\frac{\widehat{p}\left( 1 - \widehat{p} \right)}{N}}\). With these refinements the simulation matches the accelerated series across Tables 1--7, each analytic value lying within the confidence interval of its simulated counterpart, validating Proposition 1 and the associated corollaries.

The same bridge correction substantially reduces the leading time-step bias. Table 5 reports the survival probability \(P(\tau > 0.5)\) for the first example at three successively halved step sizes \(\Delta t\), with and without the correction. The uncorrected estimate carries the expected \(O\left( \sqrt{\Delta t} \right)\) overshoot, its bias roughly halving as \(\Delta t\) is quartered, while the corrected estimate is flat to within the sampling error and agrees with the analytic value (2.6). All other simulations in this section use \(\Delta t = 0.002\).

\begin{table}[htbp]
\caption{Step-size convergence of the survival probability $\mathbb{P}(\tau > 0.5)$ for the first example, computed by bridge-corrected Monte Carlo with $5 \times 10^{5}$ paths. The uncorrected estimate carries an $O(\sqrt{\Delta t})$ overshoot, while the bridge-corrected estimate is flat to within the sampling error and consistent with the analytic value (2.6). The step $\Delta t = 0.002$ is used throughout the section.}\label{tab5}
\centering
\begin{tabular}{@{}lll@{}}
\toprule
$\Delta t$ & uncorrected & bridge-corrected \\
\midrule
0.004 & 0.3469 & 0.3165 \\
0.002 & 0.3384 & 0.3172 \\
0.001 & 0.3315 & 0.3173 \\
analytic (2.6) & --- & 0.3178 \\
\botrule
\end{tabular}
\end{table}

To show the dependence on the starting point we repeat the computation from \(r_{0} = 1.55\), \(\varphi_{0} = 1.70\), nearer the outer arc and toward the \(\varphi_{2}\) edge, all other data unchanged. The survival probability again reproduces the simulation (Table 6). The boundary-hitting probabilities (Table 7) shift as anticipated: exit through the \(\varphi_{2}\) edge rises from \(\pi_{2} \approx 0.008\) in the first example to \(\pi_{2} \approx 0.16\), now a substantial and cleanly resolved fraction, while inner-arc exit remains dominant, the inward drift being unchanged. The same \emph{\ensuremath{\varepsilon}}-acceleration is applied; the accelerated probabilities agree with the bridge-corrected simulation to within its confidence interval, and their sum reproduces unit mass to \(2 \times 10^{- 4}\).

\begin{table}[htbp]
\caption{Survival probability for the second example, started nearer the outer arc at $r_{0} = 1.55$, $\varphi_{0} = 1.70$, from the series (2.6) and from $5 \times 10^{5}$ bridge-corrected trajectories (95\% confidence intervals shown).}\label{tab6}
\centering
\begin{tabular}{@{}lll@{}}
\toprule
$T$ & analytic, eq.\ (2.6) & Monte Carlo (95\% CI) \\
\midrule
0.1 & 0.9966 & $0.9966 \pm 0.0002$ \\
0.2 & 0.9229 & $0.9234 \pm 0.0007$ \\
0.3 & 0.7743 & $0.7751 \pm 0.0012$ \\
0.5 & 0.4801 & $0.4802 \pm 0.0014$ \\
0.8 & 0.2153 & $0.2155 \pm 0.0011$ \\
\botrule
\end{tabular}
\end{table}

\begin{table}[htbp]
\caption{Boundary-hitting probabilities for the second example, after $\varepsilon$-acceleration, compared with the bridge-corrected simulation ($5 \times 10^{5}$ paths). The $\varphi_{2}$-edge probability is now an order of magnitude larger than in the first example and is cleanly resolved; every analytic value lies within the simulation's 95\% confidence interval, and the accelerated sum reproduces unit mass to $2 \times 10^{-4}$.}\label{tab7}
\centering
\begin{tabular}{@{}lll@{}}
\toprule
boundary piece & analytic & Monte Carlo (95\% CI) \\
\midrule
inner arc $\pi_{in}$ & 0.8149 & $0.8150 \pm 0.0011$ \\
outer arc $\pi_{out}$ & 0.0147 & $0.0147 \pm 0.0003$ \\
edge at $\varphi_{1}$ $\pi_{1}$ & 0.0070 & $0.0071 \pm 0.0002$ \\
edge at $\varphi_{2}$ $\pi_{2}$ & 0.1638 & $0.1631 \pm 0.0010$ \\
Total & 1.0002 & 1.0000 \\
\botrule
\end{tabular}
\end{table}

\subsection{The radially symmetric limit}\label{sec:5.3}

As a final consistency check, letting \(r_{1} \rightarrow 0\) and \(\varphi_{2} - \varphi_{1} \rightarrow 2\pi\) formally degenerates the radial construction toward the disk, for which Grebenkov \citeyearpar{grebenkov2015} gives the survival probability of the radial Ornstein--Uhlenbeck process. As described in Remark 3, the two-radius determinant (2.9) then collapses: the Tricomi solution is excluded by regularity at the origin, and the eigenvalue condition reduces to the vanishing of the single Kummer function \(M\left( a,1,\alpha r_{2}^{2} \right)\) on the rotationally symmetric mode of the periodic angular problem, reproducing the classical radial result. The spurious integer roots of (2.9) disappear in the same limit, since with the second solution removed the determinant is replaced by a single entire function whose only zeros are the genuine eigenvalues. The annular-sector radial equation thus coincides with that of the disk as \(r_{1} \rightarrow 0\); the disk spectrum proper, however, belongs to the periodic angular problem and not to the Dirichlet sector, whose lowest angular order tends to \(1/2\) rather than to \(0\), the inner radius being precisely what promotes the radial problem from a single-function condition to the two-function determinant (2.9).

\subsection{A reversible correlated example}\label{sec:5.4}

To illustrate Proposition 5 and the decorrelation of Proposition 3, we take a genuinely correlated, reversible planar process. Let the diffusion and drift matrices be

\[A = \begin{pmatrix}
1.0 & 0.4 \\
0.4 & 0.6
\end{pmatrix}\quad B = \begin{pmatrix}
1.08 & 0.02 \\
0.30 & 0.36
\end{pmatrix}\quad S = \begin{pmatrix}
1.2 & - 0.3 \\
 - 0.3 & 0.8
\end{pmatrix}\]

the drift arising as \(B = AS\) with \(S\) symmetric. Then \(A^{- 1}B = S\) is symmetric, so the reversibility condition \(BA = AB^{\top}\) holds; the drift \(B\) is itself non-symmetric, so the two coordinates are genuinely coupled and the process is not a pair of independent scalar diffusions.

Diagonalizing \(B\) in the \(A^{- 1}\) inner product as in Proposition 3 produces the independent factors \(y = Mx\) with mean-reversion rates \(\lambda_{1} = 0.352\), \(\lambda_{2} = 1.088\) and volatilities \(s_{1} = 0.656\), \(s_{2} = 1.081\); one verifies that \(MBM^{- 1} = \operatorname{diag}\left( \lambda_{1},\lambda_{2} \right)\) and \(MAM^{\top} = \operatorname{diag}\left( s_{1}^{2},s_{2}^{2} \right)\) to machine precision. In the decorrelated coordinates we impose the rectangle \(y_{1} \in \lbrack - 1.5,1.8\rbrack\), \(y_{2} \in \lbrack - 1.2,1.5\rbrack\) and start the process at \(y_{0} = (0.6, - 0.4)\), equivalently \(x_{0} = ( - 0.386,0.447)\) in the original coordinates.

The one-dimensional eigenvalues and eigenfunctions entering the side-probability series (4.14) are computed exactly as in Section 5.1, from the parabolic-cylinder eigenvalue condition (4.9); here the two-point determinant carries no spurious roots, the two parabolic-cylinder solutions being entire and independent for every order. Table 8 reports the four side probabilities of Proposition 5 --- the probabilities that the process leaves the rectangle through each of its edges --- evaluated from (4.14) after the same \emph{\ensuremath{\varepsilon}}-acceleration, alongside a Monte Carlo estimate from \(1.5 \times 10^{5}\) simulated paths of the decorrelated dynamics, the two factors advanced jointly by exact Gaussian transitions with time step \(2 \times 10^{- 3}\) and each factor's boundary handled by the Brownian-bridge correction.

\begin{table}[htbp]
\caption{Side probabilities for the reversible correlated example, from the series (4.14) and from a $1.5 \times 10^{5}$-path bridge-corrected simulation of the decorrelated dynamics (95\% confidence intervals shown). The accelerated series sums to unit mass, and each analytic value lies within the simulation's confidence interval.}\label{tab8}
\centering
\begin{tabular}{@{}lll@{}}
\toprule
exit edge & analytic, eq.\ (4.14) & Monte Carlo (95\% CI) \\
\midrule
$y_{1} = h_{1}$ (factor 1, upper) & 0.0681 & $0.0674 \pm 0.0013$ \\
$y_{1} = l_{1}$ (factor 1, lower) & 0.0412 & $0.0416 \pm 0.0010$ \\
$y_{2} = h_{2}$ (factor 2, upper) & 0.2403 & $0.2403 \pm 0.0022$ \\
$y_{2} = l_{2}$ (factor 2, lower) & 0.6504 & $0.6507 \pm 0.0024$ \\
total & 1.0000 & 1.0000 \\
\botrule
\end{tabular}
\end{table}

The accelerated side probabilities agree with the simulation to within its confidence interval; the very slow first factor --- its lowest eigenvalue is \(\theta_{1,1} \approx 0.07\), so that nearly half its mass survives to \(t = 12\) --- is resolved by the acceleration without difficulty. As an independent check, the survival probability \(P(\tau > T) = S_{1}(T)\, S_{2}(T)\) of Proposition 4 factorizes across the two factors; its values \(0.9439\), \(0.8150\), \(0.6535\) at \(T = 0.2,0.5,1.0\) agree with the bridge-corrected estimates \(0.9444\), \(0.8141\), \(0.6534\) (95\% half-widths \(0.0012\), \(0.0020\), \(0.0024\)), each within the confidence interval.

\section{An application: first passage of an optically trapped particle in a sector geometry}\label{sec:6}

The results of Sections 2--3 were stated for an abstract planar Ornstein--Uhlenbeck process. We now describe a physical setting in which that process, the annular-sector geometry and, above all, the joint law of the exit time and the exit boundary all acquire a direct and measurable meaning. The domain is natural for a particle held in an isotropic trap, for which radial distance and angular position are the native, physically meaningful coordinates.

\subsection{The model: a trapped colloidal particle}\label{sec:6.1}

A micron-scale dielectric particle suspended in a viscous fluid and held in an optical trap (optical tweezers) is the paradigmatic physical realization of the planar Ornstein--Uhlenbeck process \citep{jones2015,volpe2013}. For small lateral displacements the trap exerts a linear restoring force of stiffness \(k\), and in the overdamped regime --- inertia negligible against viscous drag, as holds for a colloid in water --- the lateral position \(X_{t} \in \mathbb{R}^{2}\) obeys the Langevin equation \(\gamma\, dX_{t} = - k\, X_{t}\, dt + \sqrt{2\gamma k_{B}T}\, dW_{t}\), that is

\begin{equation}
dX_{t} = - \beta\, X_{t}\, dt + \sigma\, dW_{t},\quad\quad\beta = \frac{k}{\gamma},\quad\sigma^{2} = 2D = \frac{2k_{B}T}{\gamma},\quad\alpha = \frac{\beta}{\sigma^{2}} = \frac{k}{2k_{B}T},
\tag{6.1}
\end{equation}

where \(\gamma\) is the Stokes friction coefficient, \(D = k_{B}T/\gamma\) the Stokes--Einstein diffusion coefficient, \(T\) the temperature and \(k_{B}\) Boltzmann's constant \citep{volpe2013}. The stationary law of (6.1) is the Boltzmann distribution \(\propto e^{- k|x|^{2}/2k_{B}T} = e^{- \alpha|x|^{2}}\) in the harmonic well, which is precisely the weight \(\rho\) of Section 2 with \(\alpha = k/2k_{B}T\).

Two features make this the right physical model for the annular-sector construction, rather than a forced fit. First, a symmetric trap far from any wall is isotropic, with equal stiffness along the two lateral axes and uncorrelated thermal forcing, so the drift matrix is \(B = \beta I\) exactly --- the isotropy that Section 2 requires and that, by Remark 6, the curved geometry cannot do without. Trap asymmetry or hydrodynamic coupling would introduce the correlation of Section 4, sending the tractable domain back to a rectangle. Second, distance from the trap center and angular position are each intrinsically meaningful and can be bounded independently, so the annular sector is a natural domain and not an artificial one.

\subsection{The geometry and the four exit events}\label{sec:6.2}

Confine the particle laterally to an annular sector \(D = \{ r_{1} < r < r_{2},\ \varphi_{1} < \varphi < \varphi_{2}\}\) about the trap center. The four boundary pieces, indistinguishable in a purely radial treatment, are here four physically distinct first-passage events. The inner arc \(r = r_{1}\) represents capture on first entering the central disk it bounds, around the trap center --- a reactive site, an aperture or pore, an electrode, or a second, smaller trapped probe --- the narrow-capture target of reaction-kinetics models. The outer arc \(r = r_{2}\) represents escape past a detection radius, the particle having fluctuated far enough from the center to be registered as having left. The two radial edges \(\varphi = \varphi_{1},\varphi_{2}\) represent contact with confining walls bounding a wedge, realized either as nanofabricated channel walls or as repulsive optical line-barriers.

This is exactly the first-passage and narrow-escape setting of trapped-particle kinetics \citep{benichou2014,grebenkov2015,holcman2014,jones2015}, lifted from the disk or annulus to a domain with corners. The quantity of interest in that literature is which target is reached and when; the joint law of the exit time and the exit boundary --- the central result (3.7) and the side probabilities (3.10)--(3.12) --- answers exactly this, resolving capture at the target against outward escape against contact with either wall, together with the timing of each. The four-way decomposition is the complete outcome distribution that a single-particle-tracking experiment records.

The dimensionless problem solved in the first and second examples is scale-free: it fixes the shape of the trapped dynamics, while the absolute length and time scales are set by the physical constants through the thermal length \(\ell_{T} = \sqrt{k_{B}T/k}\) --- the root-mean-square positional fluctuation of a single coordinate in the trap --- and the relaxation time \(\beta^{- 1} = \gamma/k\). Stiffer traps compress the geometry and accelerate the dynamics; softer traps dilate and slow them; the exit statistics are identical.

\subsection{A concrete parameterization}\label{sec:6.3}

Take a silica microsphere of radius \(a = 0.5\,\mu m\) in water at \(T = 295\,\)K, for which the dynamic viscosity is \(\eta \approx 9.5 \times 10^{- 4}\, Pa \cdot s\), the Stokes friction \(\gamma = 6\pi\eta a \approx 9.0 \times 10^{- 9}\, N \cdot s \cdot m^{- 1}\), and the diffusion coefficient \(D = k_{B}T/\gamma \approx 0.45\,\mu m^{2}s^{- 1}\). In a trap of stiffness \(k = 2\, pN\,\mu m^{- 1}\) the relaxation rate is \(\beta = k/\gamma \approx 2.2 \times 10^{2}\, s^{- 1}\) (relaxation time \(\beta^{- 1} \approx 4.5\ \)ms), the trap parameter is \(\alpha = k/2k_{B}T \approx 2.5 \times 10^{2}\,\mu m^{- 2}\), and the thermal length is \(\ell_{T} \approx 45\,\)nm. Measuring length in the unit \(L = \left( \alpha_{\star}/\alpha \right)^{1/2}\), with \(\alpha_{\star} = 2.041\) the dimensionless value of the first example, gives \(L \approx 91\,\)nm and reproduces the data of both examples as the physical quantities of Table 9; the dimensionless volatility \(\sigma = \alpha_{\star}^{- 1/2} = 0.70\) is then automatic, confirming the consistency of the mapping.

\begin{table}[htbp]
\caption{Physical realization of the dimensionless parameters of the first and second examples, for a silica sphere of radius $0.5\,\mu m$ in water at $295\,$K in a trap of stiffness $k = 2\, pN\,\mu m^{-1}$ (length unit $L \approx 91\,$nm).}\label{tab9}
\centering
\begin{tabular}{@{}lll@{}}
\toprule
quantity & dimensionless & physical value \\
\midrule
mean-reversion rate $\beta$ & 1 & $2.2 \times 10^{2}\, s^{-1}$ ($\tau = 4.5\,$ms) \\
diffusion $\sigma^{2} = 2D$ & 0.49 & $0.91\,\mu m^{2}s^{-1}$ \\
trap parameter $\alpha = \beta/\sigma^{2}$ & 2.041 & $2.5 \times 10^{2}\,\mu m^{-2}$ \\
inner radius $r_{1}$ & 0.8 & $73\,$nm \\
outer radius $r_{2}$ & 2.2 & $201\,$nm \\
angular width $\varphi_{2} - \varphi_{1}$ & 2.0 & $115^{\circ}$ \\
start, the first example ($r_{0}$) & 1.30 & $119\,$nm \\
start, the second example ($r_{0}$) & 1.55 & $141\,$nm \\
\botrule
\end{tabular}
\end{table}

Thus the first and second examples describe a one-micron bead in a two-piconewton-per-micron trap escaping an annular sector roughly \(70\)--\(200\,\)nm in radius, subtending a \(115^{\circ}\) wedge, on a timescale of milliseconds --- well within the spatial and temporal resolution of modern optical tweezers. At this stiffness the wedge is set by nanofabrication; a trap an order of magnitude softer places the same geometry at the micron scale, where the angular walls may be drawn optically, with identical exit statistics.

\subsection{Reading the two examples}\label{sec:6.4}

The physical content of the two examples is now transparent. In the first example the particle starts near the trap center, and exit is overwhelmingly capture at the inner target (\(\pi_{in} \approx 0.90\)), with escape past the outer radius rare (\(\pi_{out} < 0.01\)): the equilibrium-seeking drift continually returns the particle toward the center, so the central target is reached almost surely and the outer detection radius almost never. In the second example the particle is released farther out and nearer the \(\varphi_{2}\) wall, and the probability of striking that wall rises by more than an order of magnitude (\(\pi_{2} \approx 0.16\)); this is the experimentally relevant regime of a particle released from a controlled, off-center position --- for instance after photo-uncaging or an optical kick --- whose subsequent fate is steered by the release point. In both cases the mean time to capture or escape is a few milliseconds, set by the relaxation time \(\beta^{- 1}\). The full joint law supplies the complete capture-versus-escape-versus-wall statistics and their timing, and its sensitivity to the trap stiffness, the temperature and the geometry follows by termwise differentiation of the series.

These idealizations deserve a word. The overdamped reduction is controlled by the wide separation between the momentum relaxation time --- the particle mass divided by its drag coefficient, of order a tenth of a microsecond for a micron-scale silica sphere in water --- and the millisecond exit times of Table 9, a ratio of some four orders of magnitude; inertia is therefore negligible on the timescale of the exit problem, which justifies dropping the acceleration term in the Langevin equation. The angular and radial boundaries are taken to be perfectly absorbing, an idealization of the physical capture and escape surfaces, and the drag is taken constant, neglecting the position-dependent hydrodynamic friction (Fax\'en corrections) that becomes appreciable within about one particle radius of a wall. Both approximations are accurate in the regime considered here, in which the particle explores the interior of the sector and reaches a boundary only at the terminal exit event.

\section{Conclusion}\label{sec:7}

We have obtained explicit eigenfunction expansions for the first-exit functionals of a planar Ornstein--Uhlenbeck process from two canonical domains. For the isotropic process on an annular sector, the survival probability, the exit-time density and moments, and the joint law of the exit time and the exit boundary are given by double series in elementary and confluent hypergeometric (Whittaker) functions of non-integral order, the radial spectrum being fixed by a two-radius determinant. The annular sector is the minimal domain combining curved and rectilinear boundaries; it specializes to the bounded sector and shares its radial Kummer equation with the disk, annulus and wedge. The radially symmetric survival probability of Grebenkov is recovered not by a literal Dirichlet-sector limit but by replacing the absorbing radial edges with periodic angular conditions and sending the inner radius to zero, the rotationally symmetric mode of that periodic problem reproducing the classical radial result. For a genuinely correlated, reversible process the analogous problem separates on a principal-axis rectangle, where the special functions are parabolic cylinder functions; there the survival probability factorizes into one-dimensional problems, and the non-factorizing content resides in the joint law of the exit time and exit side. The two geometries are complementary, reflecting a structural dichotomy: polar separation on the circular annular sector requires the (centered) process to be isotropic, while a reversibly correlated process separates instead on a rectilinear, principal-axis domain.

The expansions are explicit and require only one-dimensional root-finding and quadrature, the eigenvalues being isolated once the spurious roots of the determinant --- those at non-positive integer order, where the two confluent hypergeometric solutions coincide --- are discarded. We have validated the results against Monte Carlo simulation: the survival probability is reproduced throughout, and the joint exit law is reproduced to the accuracy permitted by the slower, boundary-layer-sensitive convergence of boundary functionals. Sensitivities with respect to any parameter follow by termwise differentiation, and the formulae furnish benchmarks on a cornered geometry where finite-difference and simulation methods are least reliable.

Several directions remain open. The irreversible case \(BA \neq AB^{\top}\) --- in particular a rotational, spiralling drift, for which the drift matrix has complex eigenvalues and the generator is no longer self-adjoint --- lies outside both separable frameworks and is the natural frontier. The correlated process on a curved domain is a second: it no longer separates in the polar coordinates used here, since after whitening the boundary becomes an elliptic sector, so that new special-function or numerical methods would be required except in special aligned cases. Finally, time-dependent trapping and moving boundaries, together with the optimal-stopping and control problems suggested by the statistical-arbitrage and confined-particle applications, would extend the present results in directions of evident practical interest.

\bmhead{Acknowledgments}
Not applicable.

\section*{Declarations}

\begin{itemize}
\item \textbf{Funding.} No funding was received for conducting this study.
\item \textbf{Competing interests.} The author declares no competing interests.
\item \textbf{Data availability.} No datasets were generated or analysed during the current study; all numerical results reported here are reproducible from the formulae and parameter values stated in the text.
\item \textbf{Author contributions.} T.\ Guillaume is the sole author and conducted all aspects of the work.
\end{itemize}

\bibliography{references}

\end{document}